\documentclass{amsart}
\usepackage{ulem}

\usepackage{amsmath,amssymb,amsthm,amsaddr}
 
\usepackage[absolute,overlay]{textpos}

\usepackage{epsfig}
\usepackage{epstopdf}
\renewcommand{\epsilon}{\varepsilon}
\usepackage{esint}

\usepackage{fancyhdr}
\usepackage{enumerate}
\usepackage{tikz}
\usepackage{wrapfig}
\usepackage{siunitx}
\usepackage{subfig}
\usepackage{caption}
\usepackage{booktabs}
\usepackage{braket}
\usepackage{multirow}
\usepackage{url}

\usepackage{graphicx}
\usepackage{subfig}
\usepackage{tabularx}
\usepackage{stmaryrd}

\usepackage{mathrsfs}     
\usepackage{helvet}         
\usepackage{courier}        
\usepackage{type1cm}      
\usepackage{color}
\usepackage{todonotes}
\usepackage{mathbbol}

\usepackage{geometry,calc,color}

\usepackage{enumerate}

\newcommand{\ta}{\widetilde a^K}
\newcommand{\tA}{\widetilde A_K}

\newtheorem{theorem}{Theorem}[section]
\newtheorem{corollary}[theorem]{Corollary}
\newtheorem{lemma}[theorem]{Lemma}

\theoremstyle{definition}

\newtheorem{problem}{Problem}[section]
\newtheorem{assumption}{Assumption}[section]
\newtheorem{comment}{Comment}

\theoremstyle{remark}
\newtheorem{remark}[theorem]{Remark}

\numberwithin{equation}{section}

\definecolor{ForestGreen} {cmyk}{0.91,0,0.88,0.12}

\newcommand{\rosa}[1]{{#1}} 

\newcommand{\tb}{\|}

\renewcommand{\phi}{\varphi}

\renewcommand{\hat}[1]{\widehat{#1}}

\newcommand{\shapereg}{\rho^\star}

\newcommand{\LK}{\Lambda_h}
	\newcommand{\Riz}{\mathcal{R}_K}
\newcommand{\Hstar}{H^1_{\o}}

\newcommand{\roundPrecision}{2}
\newcommand{\he}{h_e}
\newcommand{\He}{H_e}

	\newcommand{\D}{\mathcal{D}_K}

\begin{document}


\begin{textblock*}{14cm}(2.4cm,26.4cm) 
	 \noindent
\scriptsize{This is the postprint version of the paper published in {\it ESAIM: Mathematical Modelling and Numerical Analysis}, 55 (2021), S785--S810, doi: {10.1051/m2an/2020059}. The original publication is available at www.esaim-m2an.org.
	}
\end{textblock*}

	\title{{A Polygonal Discontinuous Galerkin method \\ with minus one stabilization}}\thanks{This paper has been realized in the framework of ERC Project CHANGE, which has received funding from the European Research Council (ERC) under the European Union’s Horizon 2020 research and innovation programme (grant agreement No 694515)}
	\author{Silvia Bertoluzza and Daniele Prada}
	\address{IMATI ``E. Magenes'', CNR, Pavia (Italy)}%

	%
	\begin{abstract}	We propose a Discontinuous Galerkin method for the Poisson equation on polygonal tessellations in two dimensions, stabilized by penalizing, locally in each element $K$, a residual term involving the fluxes, measured in the norm of the dual of $H^1(K)$. The scalar product corresponding to such a norm is numerically realized via the introduction of a (minimal) auxiliary space  inspired by the Virtual Element Method. Stability and optimal error estimates in the broken $H^1$ norm are proven under a weak shape regularity assumption allowing the presence of very small edges.
		The results of numerical tests confirm the theoretical estimates. 
		\end{abstract}
	%
	%
	%
	%
	\maketitle


	\maketitle
	\section{Introduction}
	
	%
	%
	%
	%
	%
	%
	%
	%
	%
	%
	%
	%
	
	Methods for solving PDEs based on polyhedral meshes are attracting more and more attention, resulting in a fast development. They provide greater flexibility in mesh generation,
	can be exploited as transitional elements in finite element meshes, and are better suited than methods based on tetrahedral or hexahedral meshes for many applications on complicated and/or moving domains \cite{DG-Rev}. Many different approaches exist, such as   the Agglomerated Finite Element method \cite{Bassietal12}, the Virtual Element Method \cite{basicVEM}, the Hybrid High Order method \cite{DiPietroErnLemaire}, just to quote the most recent. 
	
	A common ingredient to all of these methods is the presence of some stabilization term that penalizes a residual in some mesh dependent norm \cite{BCMS}. Dealing with such terms in the analysis usually relies on the use of some kind of inverse inequality, and results in suboptimal estimates when the factor stemming from such inequality does not cancel out with some small factor coming from the approximation properties of the involved space. This is the case when, for instance, the elements are not shape regular or when we want to obtain $hp$ estimates \cite{CangianietalBook,GuzmanRiviere}. This kind of problem naturally arises when a mesh dependent norm is used to mimic the action of the norm of the space where the penalized residual naturally ``lives'', usually a negative or fractionary norm. On the other hand, it has been observed that, at least theoretically, it is possible to design stabilization terms based on such a ``natural'' norm \cite{BBStab,BStab}, for which the analysis does not require the validity of any inverse inequality. 
	
{	In the following we propose a Discontinuous Galerkin method for the Poisson equation on a polygonal tessellation in two dimensions with an element by element stabilization similar to the one proposed by \cite{EwingWangWang, BurmanHansbo}, that penalizes the residual on the flux, the main novelty being the norm in which such residual is penalized, namely, the norm of the dual of $H^1$.} The numerical realization of the $(H^{1})'$ norm has been the object of several papers \cite{BPV, Arioli}, and we follow here the general approach proposed by \cite{AlgStab}. While in this paper we start by addressing the case of a mesh satisfying a weak shape regularity assumption, and we only perform the analysis of the convergence in $h$, we believe that this approach (which can, of course, be applied also to other formulations and to other problems) has the potential to tackle more general cases. 
	
	{
			The paper is organized as follows: in Section \ref{sec:2} we present and analyze the new method. More specifically, in Section \ref{sec:2.1} we define some non standard form for the norms of some Sobolev space, which make it easier to deal with the scaling of negative norms; in Section \ref{sec:2.2} we present the method, in Section \ref{sec:2.3} we define the global broken norms that we will employ in the analysis, which we carry out in Sections \ref{sec:2.4}. A separate section, namely Section \ref{proofoflemma}, is dedicated to the proof of a key inf-sup condition (Lemma \ref{lem:infsup}).
		 Section \ref{sec:3} is devoted to the definition of a computable scalar product for the dual space of $H^1$. Finally, Section \ref{sec:hybrid} presents an equivalent hybridized version of the discrete problem, particularly well suited for efficient implementation, and Section \ref{sec:4} presents the result of some numerical experiments, confirming the validity of the theoretical convergence estimate.

			As we do not aim at tracking the dependence of the constants in the estimates that we are going to provide on the polynomial degree $k$ but only on the different mesh size parameters, in order to avoid the proliferation of constants, in the following we will write $A \lesssim B$ (resp. $A \gtrsim B$) to indicate that the quantity $A$ is less or equal (resp. greater or equal) than  the quantity $B$  times a constant independent of the element diameters $h_K$, and of the edge lengths $\he$,   but possibly depending on the constant $\rho^\star$ involved in the shape regularity Assumption \ref{ass.shape.reg} and on the degree $k$ of the polynomial spaces considered.

		}

	\newcommand{\bpi}{\tilde \pi}

	\newcommand{\Tess}{\mathcal{T}_h}
	
	\newcommand{\Edges}{\mathcal{E}_h}
	\newcommand{\EdgK}{\mathcal{E}^K}
	\newcommand{\Squel}{\Sigma}
	
	\section{The DG method with minus one stabilization}\label{sec:2}
	
	\subsection{Scaled norms, seminorms and duals}\label{sec:2.1} In the following, for $\phi \in V$ and $F \in V'$ (depending on the context, $V$ and $V'$ will be different couples of dual Sobolev spaces), we will indicate by $\langle F ,\phi \rangle$ the action of $F$ on $\phi$. 
	In the analysis that follows we will rely on non standard forms for the norms of some Sobolev space. More precisely, let $D$ be a bounded Lipschitz  domain in $\mathbb{R}^d$, $d = 1,2$. We denote by $\| \cdot \|_{0,D}$ the $L^2(D)$ norm and, for $0 < s \leq 1$ we let $| \cdot |_{s,D}$ denote the $H^s(D)$ semi norm: 
	\begin{gather}\label{edgeseminorm}
\| \phi \|_{0,D}^2 = \int_D |\phi|^2, \qquad | \phi |^2_{1,D} = \int_D | \nabla \phi |^2,\\ | \phi |^2_{s,D} = \int_D \int_D \,dx\,dy \frac{| \phi(x)-\phi(y)|^2}{| x - y |^{2s+d}}, \quad 0< s < 1.
\end{gather}
Let $\sigma_D$ and $\tau_D$ be two positive constants, whose choice will be specified later. {Letting $\fint_D \phi \in \mathbb{R}$ denote the average of $\phi$ in $D$
\[
\fint_D \phi = | D |^{-1} \int_D \phi,
\]}
 we let the norm for $H^s(D)$, $0< s \leq 1$, be defined as
\[
\| \phi \|_{s,D}^2 = \sigma_D  { | \fint_D \phi |^2 }+ | \phi |_{s,D}^2.
\]
On the  dual space $(H^s(D))'$, we introduce a seminorm, defined as
	\begin{equation}\label{seminormA}
| F |_{-s,D} = \sup_{{\phi \in H^{s}(D)}\atop{\int_D \phi = 0}}\, \frac{ \langle F , \phi \rangle }{\| \phi \|_{s,D}} = \sup_{{\phi \in H^{s}(D)}\atop{\int_D \phi = 0}}\, \frac{ \langle F , \phi \rangle }{| \phi |_{s,D}},
\end{equation}
and a norm
	\begin{equation}\label{defedgedualnorm}
\| F \|^2_{-s,D} = \tau_D {| \langle F,1\rangle |^2 }+ | F |^2_{-s,D}  \end{equation}
(recall that for $s \geq 0$ the function assuming identically value $1$ over $D$ is in $H^s(D)$, so that $\langle F,1 \rangle$ is well defined).  We have the following duality result.	

\begin{lemma}\label{lem:abduality}
	Let $\sigma_D$ and $\tau_D$ satisfy {$\sigma_D \tau_D = 1$}. Then it holds
	\[
	\| F \|_{-s,D} = \sup_{\phi \in H^{s}(D)}\frac{\langle F,\phi \rangle }{\| \phi \|_{s,D}}. 
	\]
\end{lemma}

	\begin{proof} 
	Let $F\in(H^s(D))'$. We let $\bar F $ (resp. $\bar\phi$) denote, by abuse of notation, both the scalar $\bar F = \langle F,1 \rangle$ (resp. $\bar\phi = \fint_D \phi$) and the $L^2(D)$ function assuming identically 
	the value $\bar  F$ (resp. $\bar\phi$) on $D$. Observe that for all $\phi \in H^s(D)$ we have the identity
	\[
	\langle F,\phi \rangle = \bar F \bar \phi + \langle F - \bar F,\phi-\bar \phi \rangle.
	\]
Then we have
	\begin{gather*}
	\sup_{\phi \in H^{s}(D)} \frac{\langle F, \phi \rangle }{\| \phi \|_{s,D}} \leq 
	\sup_{\phi \in H^{s}(D)} 
	\frac{ |\bar F|| \bar \phi| + | F-\bar F|_{-s,D} |\phi-\bar \phi|_{s,D}}
	{\sqrt{\sigma_D | \bar \phi |^2 + |\phi|_{s,D}^2}}\leq	\\
	\sup_{\phi \in H^{s}(D)} 
	\frac
	{
		\sqrt{ \sigma_D^{-1}| \bar F |^2 + | F |^2_{-s,D}} \sqrt{\sigma_D |\bar \phi|^2 + | \phi |_{s,D}^2}
	}
	{\sqrt{\sigma_D | \bar \phi |^2 + |\phi|_{s,D}^2}} = 	{
		\sqrt{\tau_D| \bar F |^2 + | F |^2_{-s,D}}.
	}
	\end{gather*}
	On the other hand, setting $F^0 = F-| D |^{-1} \bar F$ we observe that, by the definition of $| \cdot |_{-s,D}$, for each $\epsilon >0$ there exists $\phi^0_\epsilon \in H^{s}(D)$, with $\int_D \phi_\epsilon^0 = 0$ and with $| \phi_\epsilon^0 |_{s,D} = | F^0 |_{-s,D}$, such that 
	\[
	\langle F^0, \phi^0_\epsilon \rangle \geq (1-\epsilon) | F^0 |^2_{-s,D}.
	\]
	Letting $\phi_\epsilon =  {\sigma^{-1}_D}\bar F + \phi_\epsilon^0$ we have
	\[
	\| \phi_\epsilon \|^2_{s,D} = \sigma_D^{-1}  | \bar F |^2 +  | \phi_\epsilon^0 |^2_{s,D} = \tau_D | \bar F |^2 + | F |_{-s,D}^2 = \| F \|^2_{-s,D}
	\]
	and
	\[ \langle F,\phi_\epsilon \rangle = \tau_D | \bar F |^2 + 	\langle  F^0 , \phi_\epsilon \rangle \geq (1-\epsilon) \| F \|_{-s,D}^2 = (1-\epsilon) \| F \|_{-s,D} \| \phi_\epsilon \|_{s,D}.
	\]	
	The arbitrariness of $\epsilon$ yields the thesis.
\end{proof}

\newcommand{\Nstar}{N^\star}
	\newcommand{\NV}{N_K}

	\
	
	 Let now $K \subset \mathbb{R}^2$ denote a polygon of diameter $h_K$.  More precisely, we make the following assumption, \rosa{which is quite standard in the framework of polygonal discretizations.}
	 
{	 \begin{assumption}\label{ass.shape.reg} Shape regularity: there exists a
		constant $\rho^\star >0$ such that
 $K$ is star shaped with respect to all the points in a disc of diameter $\geq \shapereg h_K$.
	 \end{assumption} }
\rosa{For the precise definition of domain star shaped with respect to a disc see \cite{Mazya:book}.} Observe that polygons for which this assumption holds satisfy \rosa{(see \cite{BrennerVEM2})}
	 \begin{equation}\label{boundperimeter}
	 | \partial K | \simeq h_K,\qquad \text{ and } \qquad | K | \simeq  h_K^2,
	 \end{equation}
	 the hidden  constants  depending on $K$ only through $\shapereg$.

	 \

{Remark that we do not make any assumption on the length of the edges of $K$ (which is not assumed to be larger than a constant times $h_K$, but is allowed to be arbitrarily small),
 or on their number (which, at least for now, we allow to be arbitrarily large), 
	 so that our assumption is weaker than what is usually assumed when dealing with the analysis of polytopal methods.
Only later on (see Section \ref{sec:3}) we will need to assume that the number of edges of the elements $K$ of the tessellation is bounded by a constant $N^\star$. 
} 
	
	\
	
	Assumption \ref{ass.shape.reg} is  sufficient to have some classical bounds with constants depending on $K$ only through $\shapereg$ 
 (see \cite{BLR,BrennerVEM2}). More precisely we have the following bounds.

\subsubsection*{Trace theorems} {For functions $u \in H^s(K)$, $1/2 < s \leq 1$ we have
\begin{equation}\label{agmon}
| u |^2_{s-1/2,\partial K} \lesssim | u |^2_{s,K}, \qquad \| u \|^2_{0,\partial K} \lesssim h_K^{-1} \| u \|^2_{0,K} + h^{2s-1}_K | u |^2_{s,K},
\end{equation}
the constant in the inequality depending on $s$. This bound is proven in \cite{BLR,BrennerVEM2} for $s = 1$, but the argument therein, based on  the
existence of a Lipschitz isomorphism $\Phi: B_1 \to K$, $B_1$ denoting the unit ball, with $\| \Phi \|_{W^{1,\infty}(B_1)} \simeq h_K$, $\| \Phi^{-1} \|_{W^{1,\infty}(K)} \simeq h_K^{-1}$, 
 applies unchanged also for $s \in (1/2,1)$, thanks to
 the 
 boundedness, for $s > 1/2$, of the trace operator from $H^s(B_1)$  to $H^{s-1/2}(\partial B_1)$.
For $u \in H^{1+s}(K)$, $1/2 < s \leq 1$ this implies that, letting $\nu_K$ denote the outer unit normal to $K$,
\begin{equation}\label{secondtraceth}
\left\| \frac{\partial u}{\partial \nu_K} \right\|^2_{0,\partial K} \lesssim  h_K^{-1} \| \nabla u \|^2_{0,K} + h^{2s-1}_K |  \nabla  u |^2_{s,K} =  h_K^{-1} |  u |^2_{1,K} + h^{2s-1}_K |  \nabla  u |^2_{s,K}.
\end{equation}}
On the other hand, for $u\in H^1(K)$ satisfying $-\Delta u = 0$ in $K$, we have that 
\begin{equation}\label{traceTh2}
| u |_{1/2,\partial K} \gtrsim | u |_{1,K}.
\end{equation}

\subsubsection*{Poincar\'e Wirtinger inequality} For $u\in H^1(K)$ we have
\begin{equation}\label{poincare}
\| u - \fint_K u \|_{0,K} \lesssim h_K | u |_{1,K}, \qquad \| u - \fint_{\partial K}  u \|_{0,K} \lesssim h_K | u |_{1,K}.
\end{equation}
	
	\

In view of Lemma \ref{lem:abduality}, on $H^s(D)$ and $(H^s(D))'$ we consider the following couple of dual norms:
\begin{equation}\label{defnorm}
\| \phi \|^2_{s,D} = {| \fint_D \phi |^2} + | \phi |^2_{s,D}, \qquad  \| F \|^2_{-s,D} = {| \langle F, 1 \rangle |^2} + | F |_{-s,D}^2.
\end{equation}

\

With these definitions, a trace theorem holds with constants only depending on the shape regularity parameter $\shapereg$.
\begin{theorem}
	It holds that 
	\[
 \| \phi \|_{1/2,\partial K} \simeq 	\inf_{{u \in H^1(K)}\atop{u=\phi \text{ on }\partial K}} \| u \|_{1,K}.
	\]
\end{theorem}

	\begin{proof}  Letting $\bar u = \fint_K u$ and $\bar u^{\partial K} = \fint_{\partial K} u$ denote the average of $u \in H^1(K)$ respectively on $K$ and on $\partial K$, we can write, thanks to \eqref{agmon},
		\[
		\| u \|^2_{1/2,\partial K} = | \bar u^{\partial K} |^2 + | u |^2_{1/2,\partial K } \lesssim | \bar u^{\partial K} |^2 + | u |^2_{1,K }.
		\]
We now observe that, as $\bar u$  coincides with the $L^2(K)$ projection of $u$ on the constants, using the boundedness of said $L^2$ projection and \eqref{poincare} we can write
\begin{equation}\label{boundbaudK}
| \bar u^{\partial K} |^2 \lesssim | \bar u |^2 + | \bar u^{\partial K} - \bar u |^2 \lesssim | \bar u |^2 + h_K^{-2}\| u -  \bar u^{\partial K} \|_{0,K}^2 \lesssim | \bar u |^2 + |  u |_{1,K}^2,
\end{equation}
which yields the first half of the thesis. As far as the second half of the thesis is concerned, letting $\phi \in H^{1/2}(\partial K)$, we let $u$ be the harmonic lifting of $\phi$.
Letting $\bar \phi = \fint_{\partial K}\phi$, and using \eqref{traceTh2} and \eqref{boundbaudK}, we have 
\[
\| u \|_{1,K}^2 = | \bar u |^2 + | u |_{1,K}^2 \lesssim | \bar \phi |^2 + | \bar u - \bar \phi |^2 + | u |_{1,K}^2 \lesssim | \bar \phi |^2 + | u |_{1,K}^2 \lesssim \| \phi \|^2_{1/2,\partial K}.
\]
which gives us the second half of the thesis.
	\end{proof}
	
\

	Remark that for $F \in (H^1(K))'$ with $\langle F,1 \rangle = 0$, the seminorm $| \cdot |_{-1,K}$ can indifferently be defined by taking the supremum over all $\phi$ with zero average on $K$ or on $\partial K$: 
	\[
	\langle F,1 \rangle = 0 \qquad \text{ implies }\qquad	\sup_{{\phi \in H^{1}(K)}\atop{\int_{\partial K}  \phi= 0}} \frac{\langle F,\phi \rangle}{| \phi |_{1,K}} = 
	\sup_{{\phi \in H^{1}(K)}\atop{\int_{K}  \phi= 0}} \frac{\langle F,\phi \rangle}{| \phi |_{1,K}}.
	\]

	Then, letting $\gamma_K: H^1(K) \to H^{1/2}(\partial K)$ denote the trace operator, and letting $\gamma^*_K$ denote its adjoint, if for $\lambda \in H^{-1/2}(\partial K)$ we have $\langle \lambda , 1 \rangle = 0$, then it holds that
	\begin{equation}\label{-1vs-1/2}
	| \gamma_K^* \lambda |_{-1,K} = \sup_{{\phi\in H^1(K)}\atop{\int_K \phi = 0}}   \frac {\langle \gamma^*_K \lambda,  \phi \rangle} { | \phi |_{1,K}} = \sup_{{\phi\in H^1(K)}\atop{\int_{\partial K} \phi = 0}}   \frac {\langle  \lambda, \gamma_K \phi \rangle} { | \phi |_{1,K}} \simeq | \lambda |_{-1/2,\partial K}, 
	\end{equation}
	(where,  $\langle \lambda , 1 \rangle$ and $\langle  \lambda, \gamma_K \phi \rangle$ stand for the duality pairing between $H^{-1/2}(\partial K)$ and $H^{1/2}(\partial K)$, while $\langle \gamma_K^* \lambda,\phi \rangle$ stands for the  duality pairing between $H^{-1}(K)$ and $H^{1}(K)$).

\

	\subsection{The model problem and its discretization}\label{sec:2.2}
Letting $\Omega$ denote a polygonal domain,	in the following we consider the simplest model problem, namely
	\begin{problem}\label{Pb:strong} Given $f \in L^2(\Omega)$ and $g\in H^{1/2}(\partial\Omega)$, find $w$ solution to
		\[
		-\Delta w = f \text{ in }\Omega, \qquad {w = g \text{ on }\partial \Omega.}
		\]	
	\end{problem}
{We assume that $g$ satisfies suitable regularity and compatibility conditions sufficient for  the existence of an $H^2(\Omega)$ function with trace equals to $g$ on $\partial \Omega$ (such assumptions are quite technical, and we refer to \cite[Theorem 2.1]{Banasiak} for more details). }
	
\

	We look for a solution to Problem \ref{Pb:strong} by a discontinuous Galerkin method on a polygonal tessellation. More precisely, let $\Tess$ denote a  tessellation of $\Omega$ into  polygons  satisfying the shape regularity Assumption \ref{ass.shape.reg}.  
	  We let  $\EdgK$ denote the set of edges of the element $K \in \Tess$, $\Edges$ denote the set of all edges of the tessellation,  and  $\Squel = \cup_{e\in\Edges} \bar e $ denote the skeleton of the decomposition. 
%

	 \
	 
{\rosa{Letting $h_K$ denote the diameter of the element $K$}, to each edge $e \in \Edges$ we associate two different mesh size parameters:  \begin{equation}\label{defHe}
	  \he  = | e |, \qquad \text{and} \qquad \He = \max_{K: e\subset \partial K} h_K,
	  \end{equation} denoting, respectively, the length of $e$, and the diameter of the largest element having $e$ as an edge. Observe that, by the definition of $\He$  \begin{equation}\label{boundrhoe}\He^{-1} \leq h^{-1}_K \ \text{ for all } e \in \EdgK, \qquad\text{and}\qquad \He \leq h_{K^+} + h_{K^-}\ \text{ for } e\subset \partial K^+ \cap \partial K^- .\end{equation} We remark that neither do we assume that, for $ e \in \EdgK$, it holds that $\he \gtrsim h_K$, nor that, for $K^+$ and $K^-$ sharing an edge, it holds that $h_{K^+} \simeq h_{K^-}$, so that our framework allows non uniform meshes with very small edges, \rosa{and adjacent elements are not constrained to have comparable diameters.}
  }

	\
	  
{On $\Squel$ we choose a unit normal $\nu $, taking care that, on $\partial\Omega$, $\nu$ points outwards.}{ We define the jump $\Lbrack u \Rbrack$ of $u = (u^K)_K \in \prod_K H^1(K)$ by setting, for all interior edges $e$ common to two elements $K^+$ and $K^{-}$,
	\begin{equation}\label{defjumpinterior}
	\Lbrack  u \Rbrack  =  u^{K^+} \nu_{K^+} + u^{K^-} \nu_{K^-}, 
	\end{equation}
	whereas, for  $e \subset \partial K \cap \partial \Omega$, we set 
	\begin{equation}\label{defjumpboundary}
	\Lbrack u \Rbrack  =  u^K \nu_K = u^K\nu.
	\end{equation}
}{Observe that the definitions \eqref{defjumpinterior} and \eqref{defjumpboundary} can be summarized in the unified expression (valid for both interior and boundary edges)
	\begin{equation}\label{defsalto}
	\Lbrack u \Rbrack|_e = \sum_{K: e \subset\partial K} u^K \nu^K.
	\end{equation}
	We underline that the cardinality of the set $\{K: e \subset \partial K\}$ is always less than or equal to two, a property that, later on, we will implicitly use at several instances.

		}

	\

\

We now let $\D: H^1(K) \to (H^1(K))'$ be defined as
\[
\langle \D u  , v \rangle = \int_K \nabla u \cdot \nabla v,
\] and, {by abuse of notation, we let $\gamma_K^*$ denote not only the adjoint of the trace operator $\gamma_K: H^{1}(K) \to H^{1/2}(\partial K)$, but also the functional $\gamma_K^* : L^2(\Sigma) \to (H^1(K))'$,  
 defined as 
\begin{equation}\label{defgammaKstar}
\langle \gamma_K^* \lambda , v \rangle = \int_{\partial K} \lambda (\nu\cdot\nu_K)  v, \qquad \text{ for all }v \in H^1(K).
\end{equation}
Observe that, if,  for some $w \in H^2(\Omega)$, $\theta  \in L^2(\Squel)$ is the single valued trace on $\Squel$ of $\nabla w\cdot \nu$, then $\gamma_K^*$ defined by \eqref{defgammaKstar} verifies
$\langle \gamma_K^* \theta,v \rangle = \langle \partial w /\partial \nu_K ,\gamma_K v \rangle$, justifying the abuse of notation.}
We have the following Lemma.
	\begin{lemma}\label{lem:2.3}  For all $\lambda \in L^2(\Squel)$ we have 
		\[
		\| \gamma_K^*\lambda \|_{-1,K} \lesssim h_K^{1/2} \| \lambda \|_{0,\partial K}.
		\]
	\end{lemma}

	\begin{proof} We have
		\begin{multline*}
		| \gamma_K^* \lambda |_{-1,K} = 
		\sup_{{\phi \in H^1(K) }\atop{\int_K \phi = 0}}\frac{\int_{\partial  K} \lambda (\nu\cdot \nu_K) \phi}{| \phi |_{1,K}}	 
		\leq   \| \lambda  (\nu\cdot \nu_K)  \|_{0,\partial K}   \sup_{{\phi \in H^1(K) }\atop{\int_K \phi = 0}}\frac{\|\phi \|_{0,\partial K}}{| \phi |_{1,K}} \\ 
		\textcolor{black}{\lesssim } \| \lambda \|_{0,\partial K}  \sup_{{\phi \in H^1(K) }\atop{\int_K \phi = 0}}  \frac{ \sqrt{h_K^{-1} \|\phi  \|_{0,K}^2 + h_K |\phi |_{1,K}^2}} {| \phi |_{1,K}}\lesssim h_K^{1/2}  \| \lambda \|_{0,\partial K},
		\end{multline*}
		where we used \eqref{agmon} and \eqref{poincare}.	
		Moreover, using a Cauchy Schwarz inequality, thanks to \eqref{boundperimeter} we can write
		\[
		| \langle \gamma_K^* \lambda, 1 \rangle | = \left| \int_{\partial K} \lambda (\nu\cdot\nu_K) \right| \leq \| \lambda \|_{0,\partial K} \| 1 \|_{0,\partial K} \lesssim  h_K^{1/2}  \| \lambda \|_{0,\partial K},
		\]
		which concludes the proof.
	\end{proof}

\

\newcommand{\klambda}{{k'}}

	We now set, for $k \geq 1$, {and $\klambda \in \{k, k-1\}$},
	\[
	V_h = \prod_{K} \mathbb{P}_k(K), \qquad {\Lambda_h = \{
	\lambda \in L^2(\Squel): \ \lambda|_e \in \mathbb{P}_\klambda(e)   \text{ for all } e \in \Edges
	\}},
	\]
	where, for any one- or two- dimensional domain $D$, $\mathbb{P}_n(D)$ denotes the space of uni- or bi- variate polynomials on $D$ of total degree less than or equal to $n$.

	\
	
	In order to define our discrete problem, we introduce, for all $K$, a bilinear form $s_K: (H^1(K))'\times (H^1(K))' \to \mathbb{R}$, satisfying the following assumption.
	\begin{assumption}\label{SK0} For all $F,G \in (H^1(K))'$ we have
		\begin{equation}\label{sK1}
		s_K(F,G) \lesssim | F |_{-1,K} | G |_{-1,K}.
		\end{equation}
		Moreover, for all $\lambda \in \LK $
		\begin{equation}\label{sK2}
		s_K(\gamma_K^* \lambda,\gamma_K^* \lambda ) \gtrsim | \gamma_K^* \lambda |_{-1,K}^2.
		\end{equation}
		
	\end{assumption}
	
	\

	We then consider the following discrete problem, where $\alpha > 0$ and $t \in \mathbb{R}$ are two parameters independent of the tessellation (and, more specifically, independent of the $h_K$'s and $h_e$'s), and where $\int_\Omega fv$ naturally  stands for $\sum_K \int_{K} f v^K$.
{ 	\begin{problem}\label{PbGlob} Find $u = (u^K)_K \in V_h$, $\lambda  \in \Lambda_h$ such that, for all $v = (v^K)_K \in V_h$, $\mu  \in  \Lambda_h$, it holds that
%
%
%
				\begin{gather}\label{PbK1}
		\sum_K	\int_K \nabla u^K \cdot \nabla v^K - 
\int_{\Squel } \lambda \Lbrack v \Rbrack \cdot \nu		
		+ t\alpha \sum_K s_K(\D u^K - \gamma_K^* \lambda ,\D v^K) =   \int_\Omega f v + t\alpha 	\sum_K	 s_K(f,\D v^K),\\
		 \label{PbK2}
	\int_{\Squel} \mu \Lbrack u \Rbrack \cdot \nu - \alpha 	\sum_K	 s_K(\D u ^K - \gamma_K^* \lambda, \gamma_K^* \mu) 	
		= \int_{\partial \Omega} g \mu - \alpha 	\sum_K	s_K(f,\gamma_K^* \mu).
		\end{gather}
	\end{problem}

	\

We easily see that Problem \ref{PbGlob} yields a consistent discretization of \eqref{Pb:strong}. Indeed, under our assumptions, the solution $w$ to Problem \ref{Pb:strong}  satisfies
$w \in H^{3/2+s}(\Omega)$ for all $s$, $0\leq s < s_0$, $s_0 >0$ depending on the geometry of $\Omega$ (see \cite[Chapter 19]{thomeeFEMbook}). This
implies $\nabla w \in H^{1/2+s}(\Omega)$ which, in turn, implies the continuity of the normal derivative across the skeleton. We can then set $\theta = \partial w /\partial \nu$, and, thanks to the trace inequality \eqref{secondtraceth} we easily see that $\theta \in L^2(\Sigma)$.
%
%
%
Multiplying the identity $-\Delta w = f$ by $v = (v^K)_K \in \prod_K H^1(K)$ and integrating by parts elementwise we  obtain 
\begin{multline}\label{weakw}
\int_{\Omega} f v =	 \sum_K \int_K \nabla w^K \cdot \nabla v^K -  \sum_K \int_{\partial K} 
\frac{\partial w^K}{\partial \nu_K} v^K =  
\sum_K  \int_K \nabla w^K \cdot \nabla v^K  - \sum_K \int_{\partial K}\theta (\nu\cdot\nu_K) v^K \\ = 
\sum_K  \int_K \nabla w^K \cdot \nabla v^K  - \int_{\Squel}\theta  \Lbrack v \Rbrack \cdot \nu. \end{multline}
Moreover we easily see that $\D w^K - \gamma_K^* \theta = f|_K$  in   $(H^1(K))'$.  It is then not difficult to check that replacing $u$ with $(w^K)_K $ ($w^K = w|_K$) and $\lambda$ with $\theta$ in \eqref{PbK1} and $\eqref{PbK2}$ yields two identities.  }

\

{	\begin{remark}
		The role of the parameter $t$ is to allow our formulation to encompass different stabilization variants in the same unified framework. While the theory presented below allows to take any $t \in \mathbb{R}$, the relevant values of $t$ are $t=0$ (for which the stabilization is, in a certain sense, minimal, as it only affects equation \eqref{PbK2}), $t=1$ (for which the stabilization term is symmetric positive semidefinite) and $t = -1$, for which we have some cancellation that can contribute to improve the inf-sup constants on which the  forthcoming analysis relies on.
		\end{remark}}
		
\begin{remark}
{For $\alpha = 0$, Problem \ref{PbGlob} is the standard hybrid formulation at the basis of the primal hybrid method \cite{RaviartThomas}, which, in \cite{EwingWangWang}, has already been  combined with a stabilization term penalizing the residual on the fluxes. The main difference between Problem \ref{PbGlob} and the method proposed in such a paper lies in the design of the stabilization term, which, in the present paper, is based on a scalar product for the space $(H^1(K))'$, whose numerical realization will be detailed later on. }Observe that the idea of measuring the residual in an $(H^1)'$ norm is not new in the context of Discontinuous Galerkin  method. In particular, it is one of the ingredient of the ultra weak formulation considered in the Discontinuous Petrov Galerkin approach  (see, for instance, \cite{DPGII}), with which the present method has certainly a number of commonalities.
	\end{remark}
	
	\subsection{Global norms and spaces}\label{sec:2.3}

	We now define the global norms 
	which we will use in our analysis. 	On $L^2(\Sigma)$ we define the norm 
\begin{gather*}
\| \lambda \|^2_{-1/2,*} = \sum_K \| \gamma_K^* \lambda \|^2_{-1,K},
\end{gather*}	
and we let  $\Lambda$  denote the  closure of $L^2(\Squel)$ with respect to such a norm. Observe that this can be identified as a closed subspace of $\prod_K (H^1(K))'$.

\

On $V = \prod_K H^1(K)$ we consider the  following  norm
\begin{equation}\label{defnorms}
\tb u \tb^2_{1,*} = \sum_K | u^K |^2_{1,K} + { \sum_{e \in \mathcal{E}_h}}\frac{ \he } {\He }| \Lbrack \bar u \Rbrack |^2,\end{equation}
where, for $u = (u^K)_K \in \prod H^1(K)$, we, once again, let $\bar u = (\bar u^K)_K$ denote the piecewise constant function defined on each $K$ as  the average  $\bar u^K = \fint_K u^K$  of $u^K$. 	
The  following Lemma states that $\tb \cdot \tb_{1,*}$  is indeed a norm on $\prod_K H^1(K)$. 

	\begin{lemma}\label{lem:poincare} For all  $u \in \prod_K H^1(K)$, letting $\bar u = (\bar u^K)_K$ denote the piecewise constant function assuming in $K$ the value $\bar u^K = \fint_K u^K$, it holds that
	\[
		\| u \|_{0,\Omega}^2 \lesssim \sum_K h^2_K | u^K |_{1,K}^2 + {\sum_{e \in \mathcal{E}_h}}\frac { \he }{\He}| \Lbrack \bar u \Rbrack |^2.
		\]
	\end{lemma}

	\begin{proof} 
		Using \eqref{poincare}, as $\bar u$ is the $L^2(\Omega)$ projection of $u$ onto $\prod_K \mathbb{P}_0(K)$,  we have
		\begin{equation}\label{1}
		\| u \|_{0,\Omega}^2 = \| u - \bar u \|_{0,\Omega}^2 + \|  \bar u \|_{0,\Omega}^2 \lesssim \sum_K h^2 _K| u^K |_{1,K}^2 +  \|  \bar u \|_{0,\Omega}^2.
		\end{equation}\label{2}
		We then only need to bound the last term on the right hand side.	Let $z$ be the solution of  
		\begin{equation}\label{eq:auxprob}
		-\Delta z =  \bar u, \text{ in }\Omega, \qquad z = 0, \text{ on }\partial \Omega. 
		\end{equation}
{Once again,		we have that $z \in H^{3/2+s}(\Omega)$ for all  $s$ with $0\leq s < s_0$, which implies the continuity of the normal derivative across the skeleton.} We can then define 
		\[
		\bar \mu^e = \fint_e \frac {\partial z} {\partial \nu}. \]
		Then,  multiplying \textcolor{black}{\eqref{eq:auxprob}} by $\bar u$ and integrating by parts element by element, we can write
		\begin{gather*}
		\|  \bar u \|^2_{0,\Omega} = - \sum_K \int_{\partial K} \frac {\partial z} {\partial \nu_K}  \bar u^K  =
		{\sum_{e \in \mathcal{E}_h}}\int_e \bar \mu^e \nu \cdot{\Lbrack}  \bar u  {\Rbrack} \leq {\sum_{e \in \mathcal{E}_h}}   \he  | \bar \mu^e | | \Lbrack \bar u \Rbrack | \\
		\lesssim 
		({\sum_{e \in \mathcal{E}_h}}  \he  \He | \bar \mu^e |^2)^{1/2}  ({\sum_{e \in \mathcal{E}_h}}\frac { \he } {\He} | \Lbrack  \bar u \Rbrack |^2)^{1/2}.
		\end{gather*}
		
		It only remains to bound the first factor in the product on the right hand side. Thanks to \eqref{boundrhoe}, we have
		\begin{multline*}
		{\sum_{e \in \mathcal{E}_h}} \he  \He | \bar \mu^e |^2 = {\sum_{e \in \mathcal{E}_h}}\He
		\int_e | \bar \mu^e |^2\leq {\sum_{e \in \mathcal{E}_h}}\He
		\int_e \left | \frac {\partial z} {\partial \nu  } \right |^2
	\leq 	{\sum_{e \in \mathcal{E}_h}}\left(\sum_{K: e\subset \partial K} h_K \right)
	\int_e \left | \frac {\partial z} {\partial \nu  } \right |^2	 	
		\\ \lesssim
		{\sum_{e \in \mathcal{E}_h}} \left( \sum_{K: e\subset \partial K } h_{K}
		\int_e \left | \frac {\partial z} {\partial \nu_{K} }\right|^2 
		\right)
		 \lesssim \sum_K h_K
		\int_{\partial K} \left | \frac {\partial z} {\partial \nu_K } \right |^2,
		\end{multline*}
where, for the last bound, we could switch the sum on $e \in \Edges$ with the sum on $K$, since the cardinality of the set $\{K: e \subset \partial K\}$ is at most two. {	Using \eqref{secondtraceth} and the regularity theory for the solution of the Poisson problem \eqref{eq:auxprob} on polygonal domains \cite[Chapter 19]{thomeeFEMbook}, we obtain (without loss of generality we can assume that, for all $K$, $h_K \lesssim 1$)
		\[
		{\sum_{e \in \mathcal{E}_h}} \he  \He | \bar \mu^e |^2 \lesssim \sum_K | z |_{1,K}^2  +   \sum_{K}| \nabla z |^{2}_{s+1/2,K} \lesssim | z |^2_{1,\Omega} + | \nabla z |^2_{s+1/2,\Omega} \lesssim
		 \|  \bar u \|^2_{0,\Omega},
		\]}	
		which yields
		\[
		\|  \bar u \|^2_{0,\Omega} \lesssim \|  \bar u \|_{0,\Omega} \left({\sum_{e \in \mathcal{E}_h}}\frac { \he } {\He} | \Lbrack  \bar u \Rbrack |^2\right)^{1/2}.
		\]
		Dividing both sides by $ \|  \bar u \|_{0,\Omega}$ and combining with \eqref{1} we get the thesis.
	\end{proof}

	Lemma \ref{lem:poincare} implies that for all $u \in \prod_K H^1(K)$, we have
	\[
	\| u \|^2_{0,\Omega} + \sum_{K} | u |_{1,K}^2 \lesssim \tb u \tb_{1,*}^2.
	\]

\newcommand{\uu}{\mathbf{u}}
\newcommand{\vv}{\mathbf{v}}

\subsection{Stability and error estimate}\label{sec:2.4}	
In order to analyze Problem \ref{PbGlob}, let us  rewrite it in  compact form:  find $\uu  = (u,\lambda) \in \mathbb{V}_h = V_h \times \Lambda_h$, such that for all $\vv = (v,\mu) \in \mathbb{V}_h$, it holds that
	\begin{gather}
	a(\uu,\vv) = F(\vv),
	\end{gather}
	with
	\begin{multline}
	a(\uu,\vv) = a(u,\lambda;v,\mu) = \sum_K \int_K \nabla u^K\cdot \nabla v^K - \int_{\Sigma} \lambda \Lbrack v \Rbrack \cdot \nu + 
	\int_{\Sigma} \mu \Lbrack u \Rbrack \cdot \nu \\ + 
	\alpha \sum_K s_K(\D u ^K - \gamma_K^* \lambda, t\D v ^K - \gamma_K^* \mu),
	\end{multline} 
	and
	\begin{equation}
F(\vv) = F(v,\mu) = \sum_K \int_K f v^K + \int_{\partial \Omega} g \mu + \alpha \sum_K s_K(f,t\D v^K - \gamma_K^* \mu ).
	\end{equation}
	
	\

	\newcommand{\constalpha}{c_\alpha}
	
It is not difficult to check that $a$ satisfies the following continuity bound: for all $u,v \in \prod_K H^1(K)$, $\lambda,\mu \in \Lambda$,
\[
a(u,\lambda;v,\mu) \lesssim \left( \sum_K \| u^K \|_{1,K}^2 + \sum_K \| \gamma_K^* \lambda \|^2_{-1,K} \right)^{1/2}   \left( \sum_K \| v^K \|_{1,K}^2 + \sum_K \| \gamma_K^*\mu \|^2_{-1,K} \right)^{1/2}.
\]
Moreover, letting
	\[
	\tb u,\lambda \tb^2_{\mathbb{V}} = \tb u \tb^2_{1,*}  + \| \lambda \|^2_{-1/2,*}
	\]
	denote the norm on $\mathbb{V} = V \times \Lambda$, we have the following lemma.
	\begin{lemma}\label{lem:infsup} 
		The bilinear form $a$ satisfies the following properties:
		\begin{enumerate}
			\item Inf-sup condition: {for all $t \in \mathbb{R}$ there exists $\alpha_0 >0$ depending on $t$ such that, for all $\alpha$, $0 < \alpha < \alpha_0$,} it holds that	\[	\inf_{(u,\lambda)\in\mathbb{V}_h}\sup_{(v,\mu)\in \mathbb{V}_h} \frac{a(u,\lambda;v,\mu)}{\tb u,\lambda \tb_{\mathbb{V}} \tb v,\mu \tb_{\mathbb{V}}} \geq \constalpha,
			\]
	{the constant $c_\alpha$ depending on $t $ and $\alpha$ but independent of the mesh size parameters $h_K$ and $h_e$.}\label{lemma27.1}
			\item Conditional continuity: for all $u,v \in \prod_K H^1(K)$, $\lambda,\mu \in L^2(\Sigma)$,  if, for all $K$, $\int_K u^K = \langle \gamma_K^* \lambda, 1 \rangle = 0$,  then we have
			\[	a(u,\lambda;v,\mu) \lesssim \tb u,\lambda \tb_{\mathbb{V}} \tb v,\mu \tb_{\mathbb{V}}\ .
			 \]\label{lemma27.2}
		\end{enumerate}
	\end{lemma}
	The proof of Lemma \ref{lem:infsup} is quite long and technical, and we postpone it to Section \ref{proofoflemma}. 
	
\

Lemma \ref{lem:infsup}(\ref{lemma27.1}) implies uniqueness of the solution to Problem \ref{PbGlob}. As such a problem is finite dimensional, uniqueness, in turn, implies existence of the solution.

\label{sec:2.5}

\
	 
	Let now $w$ be the solution of Problem \ref{Pb:strong}, and let $\theta = \partial w/\partial \nu$ denote its derivative in the normal direction $\nu$ on $\Sigma$.  We have the following lemma.
	
\rosa{	\begin{lemma}\label{lem:errest} Assume that $u \in H^1(\Omega)$ and $\theta \in L^2(\Sigma)$. Then it holds that
		\begin{equation}\label{firsterrorbound}
		\| w - u, \theta - \lambda \|^2_{\mathbb{V}} \lesssim 
		\sum_K  \inf_{v \in \mathbb{P}_k(K)} | w - v |^2_{1,K} + \sum_{e\in \Edges} \he \inf_{\mu \in \mathbb{P}_{\klambda}(e)} \| \theta - \mu \|_{0,e}^2.
		\end{equation}
		Moreover, if $u  \in H^2(\Omega)$ then we have
		\begin{equation}\label{seconderrorbound}
		\| w - u, \theta - \lambda \|^2_{\mathbb{V}} \lesssim \sum_K
		\inf_{v \in \mathbb{P}_k(K)}  \left(  | w - v |^2_{1,K} + h_K | w - v |^2_{2,K} \right).
		\end{equation}
	If, in addition, we have that $u \in H^{k+1}(\Omega)$, then it holds that
		\begin{equation}\label{finalerrorbound}
		\tb w - u, \theta - \lambda \tb^2_{\mathbb{V}} \lesssim \sum_K h_K^{2k} | w |^2_{k+1,K}.
		\end{equation}
	\end{lemma}}

\begin{proof}	
	
	Let $(w_I,\theta_I)\in \mathbb{V}_h$ be approximations to  $w$ and $\theta$ satisfying
	\begin{equation}\label{condinterpolant}
	\int_K w_I^K = \int_K w, \ \text{ for all }K\in \Tess, \qquad \int_{e} \theta_I= \int_{e} \theta \ \text{ for all }e \in \Edges. \end{equation}
	Thanks to Lemma \ref{lem:infsup}, for $u \in V_h$, $\lambda \in \Lambda_h$ solution to Problem \ref{PbGlob}, we can write
	\begin{gather*}
	\tb u - w_I,\lambda-\theta_I \tb_{\mathbb{V}} \lesssim a( u - w_I,\lambda - \theta_I;z,\zeta)
	\end{gather*}
	for some element $(z,\zeta) \in \mathbb{V}_h$ with $\tb z,\zeta \tb_{\mathbb{V}} = 1$. As observed in Section \ref{sec:2.2}, we have
	\[
	a(w,\theta;z,\zeta) = a(u,\lambda;z,\zeta),
	\]
	yielding, by Lemma \ref{lem:2.3}(\ref{lemma27.2}),
	\begin{gather*}\tb u - w_I,\lambda-\theta_I \tb_{\mathbb{V}} 
	\lesssim a( w - w_I,\theta-\theta_I;z,\zeta) \lesssim \tb w - w_I,\theta - \theta_I \tb_{\mathbb{V}},
	\end{gather*}
	and, by triangular inequality,
	\begin{equation}\label{cea}
	\tb w - u,\theta - \lambda \tb_{\mathbb V} \lesssim \tb w - w_I,\theta - \theta_I \tb_{\mathbb{V}}.
	\end{equation}
	
	\
	
	It only remains to choose suitable approximants $w_I$ and $\theta_I$ for which we can provide a bound on  the right-hand side of expression (\ref{cea}). Let then, for each $K$, $w^K_I \in \mathbb{P}_k(K)$ denote the solution to
	\[
\int_{K}(w^K_I - w) = 0, \qquad \int_{K} \nabla(w^K_I - w)\cdot \nabla z = 0, \text{ for all } z \in \mathbb{P}_k(K).
	\]

	On the other hand on each edge $e$ of $\Sigma$, let $\theta_I|_e \in \mathbb{P}_\klambda(e)$ be defined as the $L^2(e)$ projection of $\theta$. It is easy to see that the $w^K_I$'s and $\theta_I$ thus defined satisfy \eqref{condinterpolant}, so that \eqref{cea} holds. 	
			We observe that, letting $\bar w = (\bar w^K)_K$ (resp. $\bar w_I =(\bar w^K_I)_K$) with $\bar w_K = \fint_K w$ (resp. $\bar w_I^K = \fint_Kw_I^K$), thanks to \eqref{condinterpolant} we have
	\[
 \| w - w_I \|_{1,*}^2 =		\sum_K | w - w^K_I |_{1,K}^2 + { \sum_{e \in \mathcal{E}_h}}\frac{ \he } {\He} | \Lbrack \bar w - \bar w_I \Rbrack |^2 = \sum_K | w - w^K_I |_{1,K}^2. 
	\]
	Furthermore, thanks to \eqref{-1vs-1/2} and \eqref{condinterpolant}, we have
	\[
\| \theta - \theta_I \|_{-1/2,*}^2 =	\sum_K | \gamma_K^*(\theta - \theta_I) |^2_{-1,K}  \lesssim \sum_K | (\theta - \theta_I)(\nu\cdot \nu_K) |^2_{-1/2,\partial K}.
	\]

Since $\nu\cdot\nu_K$ is constant on each edge, $\theta_I (\nu\cdot\nu_K)|_e$ coincides with the $L^2(e)$ projection of $\theta(\nu\cdot \nu_K)|_e$ on $\mathbb{P}_\klambda(e)$.	By a duality argument we can then bound ($\psi_\pi \in \mathbb{P}_{\klambda}(e)$ denoting the $L^2(e)$ projection of $\psi$)
	\begin{multline*}
	| (\theta  - \theta_I)(\nu\cdot \nu_K) |_{-1/2,\partial K} = \sup_{{\psi \in H^{1/2}(\partial K)} \atop {\int_{\partial K} \psi = 0}} \frac{\int_{\partial K} ( \theta - \theta_I)(\nu\cdot \nu_K)\psi}{| \psi|_{1/2,\partial K}} 
\\	  =
	\sup_{{\psi \in H^{1/2}(\partial K)} \atop {\int_{\partial K} \psi = 0}} \frac{{ \sum_{e \in \mathcal{E}^K}} \int_e ( \theta - \theta_I)(\nu\cdot \nu_K)(\psi - \psi_\pi)}{| \psi|_{1/2,\partial K}} 
	\leq \sup_{{\psi \in H^{1/2}(\partial K)} \atop {\int_{\partial K} \psi = 0}} \frac{{ \sum_{e \in \mathcal{E}^K}} \| \theta - \theta_I \|_{0,e} \| \psi - \psi_\pi \|_{0,e}}{| \psi|_{1/2,\partial K}}\\
	\lesssim 
	\sup_{{\psi \in H^{1/2}(\partial K)} \atop {\int_{\partial K} \psi = 0}}  \frac{{ \sum_{e \in \mathcal{E}^K}} \he^{1/2} \| \theta - \theta_I\|_{0,e} | \psi |_{1/2,e}}{| \psi|_{1/2,\partial K}}
	 \lesssim \sqrt{{ \sum_{e \in \mathcal{E}^K}}  \he \| \theta - \theta_I\|^2_{0,e} },
	\end{multline*}
	where we bounded $\| \psi-\psi_\pi\|_{0,e}$ by a standard result on polynomial approximation, and used the fact that the squared piecewise $H^{1/2}$ seminorm $\sum_{e\in \EdgK} | \cdot |^2_{1/2,e}$ can be bound by the
squared	 $H^{1/2}(\partial K)$ seminorm.
	Observing that, in view of the definition of $\theta_I$ and $w_I$, we have
\[
\| \theta - \theta_I \|_{0,e} = \inf_{\mu \in \mathbb{P}_\klambda(e)} \| \theta - \mu \|_{0,e}, \quad \text{ and } \qquad | w - w^K_I |_{1,K} = \inf_{v \in \mathbb{P}_k(K) } | w - v |_{1,K},
\]
\rosa{we finally obtain \eqref{firsterrorbound} for $w \in H^1(\Omega)$ and $\theta \in L^2(\Sigma)$}.

Setting $\Lambda_h^K = \{
\mu \in L^2(\partial K): \ \mu|_e \in \mathbb{P}_{\klambda}(e) \ \text{ for all } e\in \EdgK
\}$, 
we can further bound the second term on the right hand side as follows:
\[
\sum_{e\in \Edges} \he \inf_{\mu \in \mathbb{P}_{\klambda}(e)} \| \theta - \mu \|_{0,e}^2 \leq
\sum_K \sum_{e \in \EdgK} \he \inf_{\mu \in \mathbb{P}_{\klambda}(e)} \| \theta - \mu \|_{0,e}^2 \leq
 \sum_{K} h_K \inf_{\mu \in \Lambda_h^K} \| \theta - \mu \|_{0,\partial K}^2.
\]

Then, since  for both possible choices of $\klambda$ (namely $\klambda = k$ and $\klambda = k-1$) it holds that $\nabla \mathbb{P}_k(K) \cdot\nu \subseteq \Lambda_h^K$, we can take $\mu = \partial v /\partial \nu_K$, which yields
\rosa{\begin{multline*}
\| w - u, \theta - \lambda \|^2_{\mathbb{V}} \lesssim 
\sum_K \left( 
\inf_{v \in \mathbb{P}_k(K)} | w - v |^2_{1,K}  + h_K \inf_{\mu \in \Lambda_h^K} \| \theta - \mu \|_{0,\partial K}^2
\right) \\
\lesssim \sum_K
\inf_{v \in \mathbb{P}_k(K)}  \left(  | w - v |^2_{1,K}  + h_K \left\| \theta - \frac{\partial v}{\partial \nu} \right \|_{0,\partial K}^2
\right),
\end{multline*}
and, using \eqref{secondtraceth}, if $w \in H^2(\Omega)$ we obtain \eqref{seconderrorbound}.
Assuming now that $w \in H^{k+1}(\Omega)$, standard estimates on polynomial approximation yield \eqref{finalerrorbound}}.
\end{proof}

Thanks to Lemma \ref{lem:poincare} we easily obtain a bound on the error in the standard broken $H^1$ norm. More precisely, we have the following corollary.
\begin{corollary} Assume that $u \in H^{k+1}(\Omega)$. Then we have
\[
\| w - u \|_{0,\Omega}^2 + \sum_K | w - u^K |_{1,K}^2 \lesssim \sum_K h_K^{2k} | w |^2_{k+1,K}.
\] 
\end{corollary}

\begin{remark} While \eqref{firsterrorbound} provides a sharper bound on the error, valid also when $w$ has minimal regularity, 
	if $w$ is sufficiently smooth, \eqref{seconderrorbound} has the advantage of completely decoupling the different elements, thus allowing to choose, independently in each element $K$, $\mu = \nabla v\cdot \nu$, thus bypassing the difficulty posed by the presence of possibly many small edges, and allowing for an error bound of the form \eqref{finalerrorbound} with hidden constant independent of the number of edges of the element $K$. 
\end{remark}

\rosa{	\begin{remark} The inf-sup constant $c_\alpha$ tends to $0$ linearly as $\alpha$ tends both to $0$ and to $\alpha_0(t)$.		
	 In turn, 
		 for $t$ going to $\pm \infty$, $\alpha_0(t)$ tends to $0$ as $|1+t|^{-2}$. Remark, however, that while the theoretical estimates given by Lemmas \ref{lem:infsup} and \ref{lem:errest} hold for any $t \in \mathbb{R}$, as already observed, the relevant values of $t$ are $t \in\{ -1,0,1\}$, so that, in practice,  $\alpha_0(t)$ behaves as a constant whose size depends on $\Omega$ and on the shape regularity of the tessellation. 
	On the other hand, in carrying out the proof of Lemma \ref{lem:infsup}, it can be checked (see also \cite{BPP}), that $\alpha_0(t)$ depends on the polynomial order $k$ of the method only through the possible dependence on $k$ of the implicit constants in Assumption \ref{SK0}.
	\end{remark}}

	\subsection{Proof of Lemma \ref{lem:infsup}}\label{proofoflemma}
	Let $(u,\lambda) \in \mathbb{V}_h$, and let 
	\[v = u -  \hat v\ \text{ with }\ \hat v = (\hat v^K)_K, \ \text{ where }\ \hat v^K =  \langle \gamma_K^ *\lambda,1\rangle ,\] and 
	\[\mu = \lambda + \beta \hat \mu\quad 
	\text{ with, on }e \in \Edges,\quad  \hat \mu = \He^{-1} \Lbrack \bar u \Rbrack \cdot \nu,
	\]
		where $\bar u$ denotes one more time the piecewise constant function assuming on each $K$ the value $\bar u^K = \fint_K u^K$. 
Remark that we have $\hat v \in \prod_K \mathbb{P}_0(K) \subset V_h$ as well as $\hat \mu \in \{ \mu \in L^2(\Sigma): \ \mu|_e \in\mathbb{P}_0(e) \text{ for all } e \in \Edges \}\subseteq \Lambda_h$. 

\

We can bound the $\mathbb{V}$ norm of $(\hat v,\hat \mu)$ as follows. Using \eqref{boundrhoe} and \eqref{boundperimeter} we can write 
\begin{multline}\label{boundvbar}
\| \widehat v \|_{1,*}^2 =   \sum_{e \in\Edges} \frac{\he}{\He} | \Lbrack \widehat v \Rbrack |^2 =
\sum_{e \in \mathcal{E}_h}  \frac{ \he } {\He}  \bigg| 
\sum_{K: e \subset \partial K} \langle \gamma^*_{K} \lambda,1\rangle \nu_K
\,\bigg|^2 \lesssim
\sum_{e \in \mathcal{E}_h}  \frac{ \he } {\He}  
\sum_{K: e \subset \partial K} \left|  \langle \gamma^*_{K} \lambda,1\rangle \nu_K
\right|^2
\\
\lesssim
\sum_K \frac 1{h_K} \left(
\sum_{e \in\EdgK}  \he 
\right) \left| \langle \gamma_K^* \lambda, 1\rangle  \right|^2
\lesssim
\sum_K  \left| \langle \gamma_K^* \lambda, 1\rangle \right|^2 .
\end{multline}
	Moreover, using Lemma \ref{lem:2.3} and \eqref{boundrhoe} we can write
\[	\| \gamma_K^* \hat \mu \|_{-1,K}^2 \lesssim h_K \| \widehat \mu \|^2_{0,\partial K} 
= 
h_K \sum_{e \in \EdgK} \int_e \He^{-2} | \Lbrack \bar u \Rbrack |^2 \lesssim \sum_{e \in \EdgK} \frac{\he}{\He} | \Lbrack \bar u \Rbrack |^2,
\]
which, adding over $K$ and recalling that each edge is counted at most twice, yields 
\begin{equation}\label{boundmubar}
\sum_K \| \gamma_K^* \hat \mu \|^2_{-1,K} \lesssim { \sum_{e \in \mathcal{E}_h}}\frac{ \he } {\He} | \Lbrack \bar u \Rbrack |^2.
\end{equation}

Combining \eqref{boundvbar} and \eqref{boundmubar} we obtain
\begin{equation}\label{boundvhatmuhat}
\| \hat v, \hat \mu \|_{\mathbb{V}} \lesssim \| u,\lambda \|_{\mathbb{V}}.
\end{equation}

\

	Now we have 
	\begin{multline*}
	a(u,\lambda;v,\mu) = \sum_K | u^K |_{1,K}^2  + \beta \int_{\Sigma} {\hat \mu}\,  \Lbrack u \Rbrack\cdot \nu  +  \int_{\Sigma} \lambda\, \Lbrack \hat v \Rbrack\cdot \nu +
	\alpha \sum_K s_K(\D u ^K - \gamma_K^* \lambda,t\D u^K - \gamma_K^* (\lambda+\beta{\hat \mu}))\\
	= \sum_K | u^K |_{1,K}^2 +  \beta\, I +II +\alpha\, III .
	\end{multline*}
Let us bound from below the terms $I$ through $III$.	Thanks to the definition of $\hat \mu$, as we easily see that for all $v \in \mathbb{V}$ it holds that $(\Lbrack v \Rbrack \cdot\nu) \nu =  \Lbrack v \Rbrack$, we have
	\begin{equation*}
I =  \int_{\Sigma} {\hat \mu}\,  \Lbrack u \Rbrack\cdot \nu	= { \sum_{e \in \mathcal{E}_h}}\He^{-1} \int_e \Lbrack u \Rbrack \cdot \Lbrack \bar u \Rbrack.
	\end{equation*}
By adding and subtracting $\Lbrack \bar u \Rbrack$ and using a Young inequality we can write
	\begin{multline*}
		\int_e \Lbrack u \Rbrack \cdot \Lbrack \bar u \Rbrack = \he  | \Lbrack \bar u \Rbrack |^2 +   \he  (\Lbrack \bar u^e \Rbrack - \Lbrack \bar u \Rbrack) \cdot \Lbrack \bar u \Rbrack \geq { \he }  | \Lbrack \bar u \Rbrack |^2 - { \he } | \Lbrack \bar u^e \Rbrack - \Lbrack \bar u \Rbrack | | \Lbrack \bar u \Rbrack | 
	\\
	\geq\frac12 { \he }| \Lbrack \bar u \Rbrack |^2 - \frac{1}{2} { \he } | \Lbrack \bar u^e \Rbrack - \Lbrack \bar u \Rbrack |^2,	\end{multline*}
	where, conventionally, we denote by  $\Lbrack \bar u^e \Rbrack = \fint_e \Lbrack u \Rbrack$ the average on $e$ of $\Lbrack u \Rbrack$. 
Using a Cauchy-Schwarz inequality and \eqref{defsalto} we can bound the last term as follows:
	\begin{multline*}
	| \Lbrack \bar u^e \Rbrack - \Lbrack \bar u \Rbrack |^2 = \left|  \he ^{-1} \int_e \Lbrack u - \bar u \Rbrack \right|^2 \leq  \he ^{-1}  \int_e | \Lbrack u - \bar u \Rbrack |^2
	= \he^{-1} \int_e \big| \sum_{K: e \subset \partial K} (u^K - \bar u^K) \nu^K \big|^2
	 \\ \lesssim 
\he^{-1} \sum_{K: e \subset \partial K}  \int_e | u^K - \bar u^K |^2,	\end{multline*}
	so that, using \eqref{boundrhoe}, and, once again, \eqref{defsalto} we get
	\[
	\sum_{e\in \Edges} \frac{ \he } {\He} | \Lbrack \bar u^e \Rbrack - \Lbrack \bar u \Rbrack |^2  \lesssim
	\sum_{e \in \Edges} \frac{1}{\He}  \sum_{K: e \subset \partial K}  \int_e | u^K - \bar u^K |^2
	\lesssim \sum_K h_K^{-1} \| u^K - \bar u^K \|_{0,\partial K}^2.
	\]
	Now, using \eqref{agmon} and \eqref{poincare} we have that
	\[
  \| u^K - \bar u^K \|^2_{0,\partial K} \lesssim 
	h_K^{-1} \| u^{K} - \bar u^K \|^2_{0,K} + h_K  | u^K - \bar u^K |^2_{1,K} \lesssim h_K  | u^K  |^2_{1,K},
	\]
	finally yielding, for some positive constant $c'$,
	\begin{equation}\label{term1}
	I  \geq \frac12 { \sum_{e \in \mathcal{E}_h}}\frac{ \he } {\He} | \Lbrack  \bar u \Rbrack |^2  - c' \sum_K | u^K |_{1,K}^2.
	\end{equation}
	
	We also observe that
	\begin{equation}\label{term2}
II =	\int_{\Sigma} \lambda\, \Lbrack \hat v \Rbrack \cdot \nu = \sum_K \int_{\partial K} \lambda (\nu \cdot \nu_k) \hat v^K = \sum_K |\langle \gamma_K^* \lambda,1 \rangle |^ 2.
	\end{equation}
	
\newcommand{\cstar}{c_*}
\newcommand{\Cstar}{C^*}	
Finally we can write
	\begin{multline*}
III = 
	 \sum_K  s_K( \gamma_K^* \lambda, \gamma_K^* \lambda ) + { t \sum_K s_K(\D u^K,\D u^K)} +
	\beta \sum_K  s_K( \gamma_K^* \lambda, \gamma_K^*{\hat \mu} ) 
	\\
	-  \beta \sum_K s_K (\D u ^K,\gamma^*_K{\hat \mu}) -  {(1+t)} \sum_K s_K(\D u ^K,\gamma_K^* \lambda) = IV + V + VI + VII + VIII.
	\end{multline*}
	We separately bound the five terms on the right hand side.	By Assumption \ref{SK0}, we have
	\[
IV = \sum_K  s_K( \gamma_K^* \lambda, \gamma_K^* \lambda ) \geq  c_1 \sum_K | \gamma_K^* \lambda |_{-1,K}^2.
	\]
Remarking that
	\begin{equation}\label{LimD}
	| \D u  |_{-1,K} = \sup_{{v \in H^1(K)}\atop{\int_K v = 0}} \frac{\int_K \nabla u \cdot \nabla v}{| v |_{1,K}} \leq | u |_{1,K},
	\end{equation}
{and  letting $t^- > 0$, $t^- = \max \{-t , 0 \}$,denote the negative part of $t$, thanks to  \eqref{sK1} 
we can write 
		\[
V = t \sum_K s_K(\D u^K,\D u^K) \geq - C^* t^- \sum_K | u^K |_{1,K}^2.
	\]}
	
	Using Assumption \ref{SK0}, as well as \eqref{boundmubar},  and applying a Cauchy Schwarz and a Young inequality, we also have, for some positive constant $c$,
	\begin{multline*}
| VI | \leq \beta	\sum_K | s_K( \gamma_K^* \lambda, \gamma_K^*{\hat \mu} ) | \leq \beta c\sum_K | \gamma_K^* \lambda |_{-1,K} |  \gamma^*_K{\hat \mu} |_{-1,K}  \leq \\
\beta 	c \left(\sum_K | \gamma_K^* \lambda |^2_{-1,K} \right)^{1/2}  \left(\sum_K |  \gamma_K^* \hat \mu |^2_{-1, K}\right)^{1/2} \\ 
	\leq \beta \epsilon  \sum_K | \gamma_K^* \lambda |^2_{-1,K} + \beta c_3(\epsilon) { \sum_{e \in \mathcal{E}_h}}\frac{ \he } {\He}| \Lbrack \bar u \Rbrack |^2, 
	\end{multline*}
	and, analogously,
	\begin{gather*}
| VII | \leq  \beta \sum_K |	s_K (\D u ^K,\gamma^*_K{\hat \mu}) | \leq c \beta \sum_K | u^K |_{1,K} |  \gamma^*_K{\hat \mu} |_{-1,K} \leq c_4  \beta ( \sum_K | u^K |_{1,K}^2 + { \sum_{e \in \mathcal{E}_h}}\frac { \he } {\He} | \Lbrack \bar u \Rbrack |^2),
	\end{gather*}
	whereas, 
	thanks to \eqref{LimD}, we have 
	\begin{gather*}
| VIII | \leq 
| 1+ t|  \sum_K | s_K(\D u ^K,\gamma_K^* \lambda) | \leq | 1+ t| c \sum_K | u^K |_{1,K} | \gamma_K^* \lambda |_{-1,K} \leq  \epsilon  \sum_K | \gamma_K^* \lambda |^2_{-1,K}  +  c_5(\epsilon,t) \sum_K | u^K |_{1,K}^2,
	\end{gather*}
finally yielding
\begin{equation}
III \geq (c_1 - (\beta+1)\epsilon )\sum_{K} | \gamma_K^* \lambda |_{-1,K}^2 - (C^* t^- +  c_5(\epsilon,t)) \sum_K | u^K |_{1,K}^2 
- \beta (c_3(\epsilon) + c_4) \sum_{e \in \Edges} \frac{\he}{\He} | \Lbrack \bar u \Rbrack |^2.
\end{equation}

The parameter $\epsilon $ is an arbitrary positive constant and $c_3(\epsilon)$ and  $c_5(\epsilon,t)$ are positive constants depending, respectively, on $\epsilon$, and on $\epsilon$ and $t$,  and both behaving as $\epsilon^{-1}$ as $\epsilon$ tends to $0$.
	Combining the previous bounds, we obtain	
	\begin{multline*}	
	a(u,\lambda;v,\mu) \geq 
\bigg(1 - c'\beta  - \alpha (C^* t^- +  \beta c_4 + c_5(\epsilon,t))\bigg)	\sum_K | u^K |_{1,K}^2  +    \sum_K |\langle \gamma_K^* \lambda,1 \rangle |^2 \\
+\beta \bigg(\frac 1 2 - \alpha( c_3(\epsilon) + c_4 )\bigg) 	{ \sum_{e \in \mathcal{E}_h}} \frac{ \he } {\He}| \Lbrack \bar u \Rbrack |^2  +
	\alpha  \bigg(c_1- (\beta + 1) \epsilon\bigg) \sum_K | \gamma_K^* \lambda |_{-1,K}^2.
	\end{multline*}	 
	
	\
	
	We now set $\beta = 1/(2c')$, and we choose $\epsilon$  in such a way that $(\beta + 1) \epsilon = c_1/2$. With this choice, it is not difficult to see that, setting 
	\[
	\alpha_0(t) =\frac 1 2 \min\left\{
 \big( 
	c_3(\epsilon) + c_4
	\big)^{-1},\left(
	\frac{c_4 }{2c'}+ c_5(\epsilon,t) + C^* t^-
	\right)^{-1}
\right\},
	\]
if $\alpha < \alpha_0(t)$, then
	\begin{equation*}	
a(u,\lambda;v,\mu) \gtrsim
\| u,\lambda \|^2_{\mathbb{V}},
\end{equation*}	 
	the implicit constant in the inequality depending on $\alpha$ and $t$. 
		Observe that neither $\beta$ nor $\alpha$ depend on the mesh size parameters $h_K$ and $h_e$. 
	Using \eqref{boundvhatmuhat}, we then get that 
	\begin{gather*}
	\sup_{(v,\mu)\in\mathbb{V}_h} \frac{a(u,\lambda;v,\mu)} {\tb v,\mu \tb_\mathbb{V}}  \geq \frac{a(u,\lambda; u - \kappa \hat{v},\lambda+\beta \hat \mu)}{\tb u - \kappa \hat{v},\lambda+\beta \hat \mu \tb_\mathbb{V}}
	\gtrsim \frac{\tb u, \lambda \tb_\mathbb{V}^2 }{\tb u, \lambda \tb_\mathbb{V}},
	\end{gather*}
	which concludes the proof of point (1).

	\
	
	Let us now consider the continuity bound (point (2)). We observe that if $\langle \gamma_K^* \lambda, 1 \rangle = 0$ we have
	\[
	\int_{\Sigma} \lambda \Lbrack v \Rbrack \cdot \nu= \sum_K \langle \gamma_K^* \lambda , v^K  \rangle = \sum_K \langle \gamma_K^* \lambda , v^K - \fint_K v^K \rangle \leq \sum_K | \gamma_K^* \lambda |_{-1,K} | v^K |_{1,K},
	\]
	while, if $\fint_K u = 0$ we can write, for all $\mu \in L^2(\Sigma)$
	\[
		\int_{\Sigma} \mu  \Lbrack u \Rbrack \cdot \nu= \sum_K \langle \gamma_K^* \mu , u^K  \rangle \leq
	\sum_K	| \gamma_K^* \mu |_{-1,K} | u^K |_{1,K} \lesssim 
	\sum_K	| \gamma_K^* \mu |_{-1,K} | u^K |_{1,K}.
	\]
Thanks to these inequalities, in view of Assumption \ref{SK0}, the continuity bound of point (2) is easily proven, by a Cauchy-Schwarz inequality.
	
	\begin{remark}
		It is not difficult to realize that the inf-sup bound holds for all subspace $\mathbb{V}_h = V_h \times \Lambda_h \subseteq \mathbb{V}$, provided 
		$V_h \supseteq \prod_K\mathbb{P}_0(K)$ and $\Lambda_h \supseteq \{\mu \in L^2(\Sigma): \mu|_e \in  \mathbb{P}_0(K)\ \forall e \in \Edges \}$.
	\end{remark}

\section{Realizing a computable stabilizing term}\label{sec:3} In order for the proposed method to be practically feasible, we need to construct a computable bilinear form $s_K$ satisfying (\ref{sK1}) and (\ref{sK2}). The numerical realization of scalar products for negative Sobolev spaces has been the object of several papers \cite{BPV,BCT00,Arioli}.
	In particular, following the approach of \cite{AlgStab}, we introduce an auxiliary space {$Y_K \subseteq \Hstar(K) = \{v \in H^1(K), \ \fint_K v = 0\}$}, satisfying
	\begin{equation}\label{mangiacoda}
	\inf_{{\lambda \in \Lambda_h^K }}\ \sup_{{y \in Y_K}} \frac{\int_K \lambda y}{| \gamma^*_K\lambda |_{-1,K} | y |_{1,K}} \gtrsim 1, \qquad\text{where }\quad \Lambda_h^K = \{\mu \in L^2(\partial K):\ \mu|_e \in \mathbb{P}_{\klambda}(e),\ \forall e\in \EdgK  \}.
	\end{equation}
	We let $\phi_i$, $i=1,\cdots,N$ denote a basis for $Y_K$, and we let  $\ta: Y_K \times Y_K \to \mathbb{R}$ denote a continuous, symmetric bilinear form satisfying, for all  $x,y \in Y_K$,
	\begin{equation}\label{propta}
	\ta(x,y) \lesssim | x |_{1,K} | y |_{1,K}, \qquad \text{ and }\qquad \ta(x,x) \gtrsim | x |_{1,K}^2.
	\end{equation}
Letting	$\tA$ denote  the corresponding stiffness matrix
	\[
	\tA = (\ta_{ij}), \text{ with } \ta_{ij} = \ta(\phi_j,\phi_i),
	\] which, thanks to the Poincar\'e inequality, is invertible,
we can now introduce the  bilinear form $s_K: (H^{1}(K))' \times (H^{1}(K))' \to \mathbb{R}$ defined as follows:
	\begin{equation}\label{defskdiscrete}
	s_K(F,G) = \vec f^T \tA^{-1} \vec g, \quad \text{ with } \vec f = (\langle F,\phi_i\rangle)_{i=1}^N, \  \vec g = (\langle G,\phi_i\rangle)_{i=1}^N.
	\end{equation}

		\
		
		Observe that for $u,v \in H^{1}(K)$, $\lambda,\mu \in L^2(\Sigma)$ and $f \in L^2(K)$ we have 
		\[
		s_K(\D u - \gamma_K^* \lambda,t \D v - \gamma_K^* \mu) = \vec \eta^T \tA^{-1} \vec \zeta, \qquad s_K (f,t \D v - \gamma_K^* \mu ) = \vec f^T \tA^{-1} \vec \zeta,
		\]
		with $\vec \eta = (\eta_i)_{i=1}^N$, $\vec \zeta = (\zeta_i)_{i=1}^N$, $\vec f = (f_i)_{i=1}^N$ given by
		\begin{equation}\label{vecetazeta}
		\eta_i = \int_K \nabla u \cdot \nabla \phi_i - \int_{\partial K} \lambda(\nu\cdot\nu_K) \phi_i, \qquad \zeta_i = t \int_K \nabla v \cdot \nabla \phi_i - \int_{\partial K} \mu(\nu\cdot\nu_K)  \phi_i, \qquad f_i = \int_K f \phi_i.
		\end{equation}
		
		\
	
{The bilinear form $s_K$ satisfies \eqref{sK1}. Indeed, since $s_K$ is symmetric positive definite,  we have a Cauchy Schwarz inequality:
\begin{equation}\label{C-SsK}
s_K(F,G) \lesssim \sqrt{s_K(F,F)} \sqrt{s_K(G,G)}.
\end{equation}
Now, given $F \in (\Hstar(K))'$, and letting $x^F = \sum_{i=1}^N x^F_i \phi_i \in Y_K$ be the solution to
\[ \ta (x^F,y) = \langle F, y \rangle,\ \text{ for all } y \in Y_K,\]
a standard argument yields, 
\[
| x^F |^2_{1,K}	\lesssim \ta(x^F,x^F)	=  \langle F,x^F \rangle \lesssim | F |_{-1,K} | x^F |_{1,K}.
\]
Dividing both sides by $| x^F |_{1,K}$ we obtain that $| x^F |_{1,K} \lesssim | F |_{-1,K}$.
We now observe that, letting $\vec x^F=(x^F_i)_{i=1}^N$ denote the vector of coefficients of $x^F$ (which is easily seen to satisfy the identity $\vec x^F = \tA^{-1} \vec f$, with $\vec f$ given by \eqref{defskdiscrete})  we have
\[
s_K(F,F) = \vec f^T \tA^{-1} \vec f = \vec f^T \vec x^F = \sum_{i=0}^N \langle F,\phi_i \rangle x^F_i = \langle F , x^F \rangle  \lesssim | F |_{-1,K} | x^F |_{1,K} \lesssim | F |_{-1,K}^2.
\]
A similar bound holds for $G$, which, combined with \eqref{C-SsK} yields \eqref{sK1}.

\

On the other hand, let $\lambda \in \Lambda_h$, and let now $x^\lambda = \sum_{i=1}^N x^\lambda_i \phi_i \in Y_K$ denote the solution to \[\ta(x^\lambda,y) = \langle \gamma_K^* \lambda, y \rangle, \ \text{ for all }
 y \in Y_K.\] Assuming that \eqref{mangiacoda} holds, and using \eqref{propta}, we can write, for some element $y^\lambda \in Y_K$,
\begin{equation}\label{step2181}
| \gamma_K^* \lambda |_{-1,K} \lesssim
\frac{\langle \gamma_K^* \lambda , y^\lambda \rangle}{| y^\lambda |_{1,K}} =
 \frac{a(x^\lambda,y^\lambda)}{| y^\lambda |_{1,K}} \lesssim \frac{\displaystyle{\sqrt{\ta(x^\lambda,x^\lambda)}\sqrt{\ta(y_\lambda,y_\lambda)}} }{| y^\lambda |_{1,K}} \lesssim \sqrt{\ta(x^\lambda,x^\lambda)}.
\end{equation}
It is now easy to check that, setting $\vec x^\lambda = (x^\lambda_i)_{i=1}^N$ and  $\vec \lambda = (\lambda_i)_{i=1}^N$, with $\lambda_i = \langle \gamma_K^* \lambda, \phi_i \rangle$, we have that $\vec x^\lambda = \tA^{-1} \vec \lambda$ and
\begin{equation}\label{step2182}
\ta(x^\lambda,x^\lambda) = (\tA^{-1} \vec\lambda)^T  \tA (\tA^{-1} \vec \lambda) = \vec\lambda^T \tA ^{-1} \vec \lambda = s_K(\gamma_K^* \lambda,\gamma_K^* \lambda).
\end{equation}
Combining \eqref{step2181} and \eqref{step2182} we easily obtain \eqref{sK2}  (actually, a stronger result holds, namely, under our assumptions on $\widetilde a$, it is possible to prove (see \cite{AlgStab}) that (\ref{mangiacoda}) is a necessary and sufficient condition for  \eqref{sK2} to hold).}

	\
	
	We then only need to choose a (small) space $Y_K$ satisfying \eqref{mangiacoda} (remark that $Y_K$ is not required to satisfy any approximation property). 
	We choose a suitable subspace of the local non conforming Virtual Element space of order $\klambda+1$ (see \cite{noncVEM}). More precisely we set
	\[
	Y_K = \{
	y \in \Hstar(K): \frac{\partial y}{\partial \nu_K}|_e \in \mathbb{P}_\klambda(e),\ -\Delta y \in \mathbb{P}_{\klambda-1}(K), \ \int_K y p = 0\ \forall p \in \mathbb{P}_{\klambda-1}(K)
	\}
	\]
{ In order to be able to work with average free functions also for $\klambda = 0$, we use the convention that $\mathbb{P}_{-1}(K) = \mathbb{P}_0(K)$, that is, we consider what, in the virtual element framework, is referred to as an {\it enhanced space}.}

\

It is not difficult to	check that $Y_K$ satisfies condition \eqref{mangiacoda}. This is a consequence of the following Lemma. 
\begin{lemma}
	For all $y \in Y_K$ it holds 
	\[
 | y |_{1,K} \lesssim	\left| \gamma_K^* \left(\frac{\partial y}{\partial \nu_K}\right) \right|_{-1,K} \lesssim  | y |_{1,K}.
	\]
	\end{lemma}
	
	\begin{proof}
Let $y \in Y_K$. Integrating by part and using the definition of $Y_K$ we have
		\begin{gather*}
| y |_{1,K}^2 = \int_K | \nabla y |^2 = - \int_K y \Delta y + \int_{\partial K} \frac{\partial y}{\partial \nu_K} y =  \int_{\partial K} \frac{\partial y}{\partial \nu_K} y \lesssim | \gamma_K^*  
\left(
\frac{\partial y}{\partial \nu_K} \right) |_{-1,K} | y |_{1,K}. 
			\end{gather*}
			Dividing both sides by $| y |_{1,K}$ we get the first of the two bounds. On the other hand we have
			\begin{gather*}
			\left| \gamma_K^* \left(\frac{\partial y}{\partial \nu_K}\right) \right|_{-1,K}  = \sup_{{\phi \in H^1(K)}\atop{\int_K \phi = 0}} \frac{\int_{\partial K} \frac{\partial y}{\partial \nu_K} \phi}{| \phi|_{1,K}} =  \sup_{{\phi \in H^1(K)}\atop{\int_K \phi = 0}} \frac{
			\int_K \Delta y  \phi + \int_K \nabla y \cdot \nabla \phi
				}{| \phi |_{1,K}}\\  \lesssim  \sup_{{\phi \in H^1(K)}\atop{\int_K \phi = 0}} \frac{\| \Delta y \|_{0,K} \| \phi \|_{0,K} + | y |_{1,K} | \phi |_{1,K}}{| \phi |_{1,K}} \lesssim | y |_{1,K},
			\end{gather*}
			where we used a Poincar\'e Wirtinger inequality \eqref{poincare}, and  an inverse inequality of the form $\| \Delta y \|_{0,K} \lesssim h_K^{-1} | y |_{1,K}$ which holds for all functions such that $\Delta y \in \mathbb{P}_{\klambda-1}(K)$, provided Assumption \ref{ass.shape.reg} holds (see \cite{BLR} for a proof).			
		\end{proof}

	In view of the previous Lemma, the inf-sup condition \eqref{mangiacoda} is then easily proven. Indeed, given $\lambda \in \Lambda^K_h$, we let $y^ \lambda \in \Hstar(K)$ denote the (unique, as the problem is well posed as shown in \cite{noncVEM}) function with
	\[
	-\Delta y^\lambda \in \mathbb{P}_{\klambda-1}(K), \qquad \frac{\partial y^\lambda}{\partial \nu_K} = \lambda \text{ on }\partial K, \qquad \int_{K} y^\lambda p = 0 \text{ for all }p \in  \mathbb{P}_{\klambda-1}(K).
	\] 
	We have $y^\lambda \in Y_K$ and
	\[
\sup_{y \in Y^K} \frac{\int_{\partial K} \lambda y}{| y |_{1,K}} \geq\frac{\int_{\partial K} \lambda y^\lambda}{| y^\lambda |_{1,K}} =
\frac{\int_{\partial K} \frac{\partial y^\lambda}{\partial \nu_K} y^\lambda}{| y^\lambda |_{1,K}}=
 \frac{\int_K y^\lambda \Delta y^\lambda  + \int_K | \nabla y^\lambda |^2}{| y^\lambda |_{1,K}}  
= | y^\lambda |_{1,K}
\gtrsim  | \gamma_K^* \lambda |_{-1,K},
	\]
	where we used the fact that, by the definition of $Y_K$, $y^\lambda$ is $L^2(K)$ orthogonal to $\Delta y^\lambda$, as the latter is a polynomial in $\mathbb{P}_{\klambda-1}(K)$.

	\
	
	\newcommand{\bli}{e_i}
	\newcommand{\blj}{e_j}
	\newcommand{\dofi}{c_i}

	For $\klambda \geq 1$, a function $y \in Y_K$ is uniquely determined by the value of its moments up to order $k'$ on each edge.  In fact, the remaining degrees of freedom for the full non conforming VEM space of order $k'+1$ are the interior moments up to order $k'-1$ (see again \cite{noncVEM}), which we fixed to be zero in the definition of $Y_K$. {Moreover, using the same arguments as in \cite[Lemma 3.1]{noncVEM} it is easy to see that, also for $\klambda=0$, a function in $Y_K$ is uniquely determined by the value of its zero order moments on each edge.} In both cases, equivalently,  a function $y \in Y_K$ is uniquely determined by the value of the $L^2(\partial K)$ scalar products with the elements of a basis $\{
	\bli, \ i = 1,\cdots,(k+1) \NV\}$ of the space $\LK^K$ ($\NV$ denotes here the number of edges of $K$). 
 
 \
 
We let $\phi_i$ denote  the unique function in $Y_K$ for which, for all $j \in\{ 1,\cdots,(k+1) \NV\}$, $\int_{\partial K} \phi_i \blj = \delta_{ij}$, so that a function $y \in Y_K$ can be expressed as 
\[
y = \sum_{i=1}^{(k+1)\NV} \dofi \phi_i \quad \text{ with } \quad \dofi = 
\int_{\partial K} y \bli. 
\]
 As customary in the Virtual Element framework, the basis functions $\phi_i$ are not explicitly known, but the knowledge of the degrees of freedom $\dofi$, $i = 1,\cdots,(k+1)\NV$  is sufficient to compute the vectors $\vec \eta$ and $\vec \zeta$. In fact, for $u \in \mathbb{P}_k(K)$ and $\lambda \in \LK $ we have
	\[
	\eta_i = \int_K \nabla u \cdot \nabla \phi_i - \int_{\partial K} \lambda (\nu\cdot\nu_K) \phi_i = - \int_K \Delta u \phi_i + \int_{\partial K}\left( 
	\frac {\partial u}{\partial \nu_K} - \lambda (\nu\cdot\nu_K) 
	\right)  \phi_i = \int_{\partial K}\left( 
	\frac {\partial u}{\partial \nu_K} - \lambda (\nu\cdot\nu_K) 
	\right)  \phi_i,
	\] 
	where we once again used that $\phi_i$ is orthogonal to all polynomials in $\mathbb{P}_{\klambda-1} (K)\supseteq \mathbb{P}_{k-2}(K)$, and hence to $\Delta u$. As both $\lambda(\nu\cdot\nu_K) $ and $\partial u / \partial \nu_K$ belong to $\LK^K$, it is possible to write them as a linear combination of the basis functions $\bli$:
	\[
	\lambda(\nu\cdot\nu_K)  = \sum_i x_i \bli, \qquad \frac{\partial u}{\partial \nu_K} = \sum_i y_i e_i.
	\]
	Then
	\[
	\eta_i = \int_{\partial K}\left( 
\frac {\partial u}{\partial \nu_K} - \lambda (\nu\cdot\nu_K)
\right)  \phi_i = y_i - x_i.
	\]
	Moreover, the fact that $\phi_i$ is orthogonal to polynomials   in $\mathbb{P}_{\klambda-1}(K)$ also allows us to approximate $f_i \approx 0$ (which corresponds to approximating $f$ in $K$ with a polynomial in $\mathbb{P}_{\klambda-1}(K)$).

 \
 
 \newcommand{\SK}{\sigma^K}
We choose $\ta$ as the non conforming Virtual Element approximation of the bilinear form $\int_K \nabla y \cdot \nabla x$. More precisely, letting $\Pi^\nabla_K: H^1(K) \to \mathbb{P}_{k'+1}(K)$ denote the projection operator defined by the conditions
 \[
 \int_K \nabla(\Pi_K^\nabla y) \cdot \nabla q =  \int_K \nabla y \cdot \nabla q, \ \forall q \in \mathbb{P}_{k'+1}(K), \qquad \text{ and }\qquad  \int_K \Pi_K^\nabla y = 0,
 \]
 we set
 \begin{equation}\label{defta}
 \ta(x,y) = \int_K \nabla \Pi^\nabla_K x \cdot\nabla \Pi^\nabla_K y + \SK(x-\Pi^\nabla_K x,y-\Pi^ \nabla_K y),
 \end{equation}
 where, for all $x$ with $\Pi^\nabla_K x = 0$, the bilinear form $\SK$ satisfies
 \begin{equation}\label{condSKVEM}
 \SK(x,x) \simeq | x |_{1,K}^2, \qquad \SK(x,y) \lesssim | x |_{1,K} | y |_{1,K}.
 \end{equation}

We recall (see \cite{noncVEM}) that $\Pi_K^\nabla y$ is directly computable for all $y \in Y^K$ as a function of the degrees of freedom $c_i$.
We are then left with the problem of choosing a computable bilinear form $\SK$ satisfying \eqref{condSKVEM}.  A possible choice for $\SK$ is the following
 \begin{equation}\label{stabVEM}
 \SK(x,y) = \sum_{e \in \mathcal{E}^K}  \he ^{-1} \int_{e} \Pi^{ \partial}_K (\gamma_K x) \Pi^{ \partial}_K (\gamma_K y),
 \end{equation}
 where $ \Pi^{ \partial}_K: L^2(\partial K) \to \Lambda_h^K$ denotes the $L^2(\partial K)$ orthogonal projection. For such a bilinear form, condition \eqref{condSKVEM} is proven in  \cite{MascPer} for all $y$ with 
 \[
 \frac {\partial y}{\partial \nu_K}  \in\Lambda^K_h, \qquad -\Delta y = 0,
 \]
 under a stronger shape regularity assumption, namely that $ \he  \simeq h_K$ for all $e \in \mathcal{E}^K$. 
  {In the more general case that we consider here, condition \eqref{condSKVEM} holds with constants only weakly depending on the ratio $h_K/h_e$,  provided that the tessellation satisfies the following additional shape regularity assumption.}
 
{ \begin{assumption}\label{Ass3.1} There exists a constant $N^\star$ such that all the elements of the tessellations $\Tess$ have at most $N^\star$ edges.
 	\end{assumption}
 
{Under such an assumption,  see \cite{BPP_NcVEM}, it can be proven that
  \begin{equation*}{
 \SK(x,x) \gtrsim \log\left(\frac{h_K}{\min_{e\subseteq \partial K} h_e}
 \right)^{-1} | x |_{1,K}^2, \qquad \SK(x,y) \lesssim | x |_{1,K} | y |_{1,K}.}
 \end{equation*}}
 
 Of course, Assumption \ref{Ass3.1} is always satisfied with $N^\star = \max_{K} \NV$, however, for $\widetilde a$ defined by \eqref{defta}, with $\SK$  given by \eqref{stabVEM}, such a value will affect the constants in \eqref{propta}.
}

{	\begin{remark} A necessary condition for an inf-sup bound of the form \eqref{mangiacoda} to hold  is that $\dim{Y_K} \geq \dim{\Lambda_h^K}$.  As in our case the dimension of $Y_K$ verifies $\dim{Y_K} = \dim{\Lambda_h^K}$, such a space is   of the minimal dimension needed for such a condition to hold. Of course, other choices are possible for the space $Y_K$. A possibility is to choose a space of supremizers (see \cite{Rozza_supremizer}), which, in our case, would be
		\[
		Y_K = \{
		y \in \Hstar(K): \frac{\partial y}{\partial \nu_K}|_e \in \mathbb{P}_\klambda(e),\ -\Delta y =0\}.
		\]
		This is also a subspace of the local non conforming VEM space of order $\klambda + 1$, so that one can build the corresponding bilinear form $s^K$ starting from the same bilinear form $\widetilde a^K$  defined by \eqref{defta}. However, though such a choice would also lead to a bilinear form $s_K$ satisfying \eqref{mangiacoda}, it would not be possible to compute the contribution of the right hand side to the stabilization term, as, for such a choice, we do not have access to the values of the moments of the basis functions $\phi_i$. \rosa{Another possible choice, which however leads to a larger auxiliary space $Y_K$, is to resort to a finite element space of order $k$ on a sufficiently fine sub-triangulation of the polygon $K$, as it is done, though in a different spirit, in \cite{Barranchea}}.
		 \end{remark}}

	\newcommand{\PiW}{\widehat\Pi^K}

\newcommand{\LKloc}{\Lambda_h^ K}
\newcommand{\hlambda}{\widehat\lambda}
\newcommand{\hlambdaK}{\widehat \lambda^K}
\newcommand{\hmu}{\widehat\mu}
\newcommand{\hmuK}{\widehat\mu^K}
\newcommand{\hF}{\widehat F}

 \section{A hybridized version of the method}\label{sec:hybrid}
{By introducing an independent approximation of the trace of $w$ on $\Sigma$, and by replacing the single valued  approximation $\lambda$ of $\partial w/\partial \nu$ with independent approximations $\hlambdaK$  of $\partial w / \partial \nu_k =  \lambda (\nu \cdot \nu_K)$, we obtain an equivalent formulation of our problem which is better suited for efficient implementation. More precisely, we set
	\[
\widehat\Lambda_h = \prod_{K} \LKloc , \text{ with }\LKloc  = \{\hlambda \in L^2(\partial K): \hlambda|_e \in \mathbb{P}_\klambda(e)  \text{ for all } e \in \EdgK\},
	\]
	as well as 
	\begin{gather}
	\Phi_{h} = \{ \phi \in L^2(\Sigma): \phi|_e \in \mathbb{P}_\klambda(e) \text{ for all } e \in \Edges,  \ \phi|_{\partial\Omega} = 0\}.
	\end{gather}
	We then consider the following discrete problem.}
	
	\
		\begin{problem}\label{PbGlobHybrid} Find $u = (u^K)_K \in V_h$, $\hlambda = (\hlambdaK)_K \in \widehat{\Lambda}_h$, $\phi \in \Phi_{h}$ such that, for all $K \in \Tess$, for all $v^K \in \mathbb{P}_k(K)$, $ \hmuK \in \LKloc$, it holds that
		\begin{gather}\label{PbKh1}
	\int_K \nabla u^K \cdot \nabla v^K - 	 \int_{\partial K} \hlambdaK v^K + t\alpha s_K(\D u^K - \gamma_K^* \hlambdaK,\D v^K) = 		 \int_K f v^K + t\alpha s_K(f,\D v^K),\\ \label{PbKh2}
			\int_{\partial K} u^K \hmuK - \alpha  s_K(\D u ^K - \gamma_K^* \hlambdaK, \gamma_K^* \hmuK) - 	\int_{\partial K} \phi \hmuK  =\int_{\partial K\cap\partial \Omega}g \widehat{\mu}^K- \alpha 	s_K(f,\gamma_K^* \hmuK), 
		\end{gather}
		and for all $\psi \in \Phi_h$
		\begin{equation}
	 \label{Pbcoupling}
		\sum_K \int_{\partial K} \hlambdaK \psi = 0.
		\end{equation}		
	\end{problem}
	
	\
	
	Observe that, for each $K$, (\ref{PbKh1}-\ref{PbKh2}) yield a local discrete 
	 problem with Dirichlet boundary conditions imposed by Lagrange multipliers, and with a non standard stabilization term, while (\ref{Pbcoupling}) imposes continuity of the fluxes $\hlambda$ across $\Sigma$.   The coupling between the different local problems stems from the common Dirichlet data $\phi$, which is single valued on the interface, as well as from equation (\ref{Pbcoupling}). 
	
	\

{Problem \ref{PbGlobHybrid} is well posed. Indeed, letting $\widehat{\mathbb{V}}_h = V_h \times \widehat \Lambda_h$, and setting
	\[
	\widehat a(u,\hlambda;v,\hmu) = \sum_K \widehat a_K(u^K,\hlambdaK;v^K,\hmuK), \quad b(v,\hmu;\phi) = \sum_K b_K(v^K,\hmuK;\phi),\quad \hat F (v,\hmu) = \sum_K \hat F_K(v^K,\hmuK),
	\]
	where the local bilinear forms $\widehat a_K$ and $b_K$, and the linear operator $\widehat F_K$ are  respectively defined as 
	\begin{multline}
	\widehat a_K(u^K,\hlambdaK;v^K,\hmuK) =  \int_K \nabla u^K\cdot \nabla v^K - \int_{\partial K} \hlambdaK v^K \\ + 
	\int_{\partial K} u^K \hmuK  + 
	\alpha s_K(\D u ^K - \gamma_K^* \hlambdaK, t\D v ^K - \gamma_K^* \hmuK),\end{multline} 
	and
	\begin{equation}
	b_K(v^K,\hmuK;\phi) = \int_{\partial K} \hmuK \phi, \quad \hat F_K(v^K,\hmuK) = \int_K f v^K + \int_{\partial K\cap \partial\Omega} g\hmuK + \alpha s_K(f,t\D v^K -\gamma_K^* \hmuK),
	\end{equation}
	Problem \ref{PbGlobHybrid}  rewrites as: find $(u,\widehat \lambda) \in \widehat{\mathbb{V}}_h$, $\phi \in \Phi_h$ such that for all $(v,\widehat\mu) \in \widehat{\mathbb{V}}_h$, $\psi \in \Phi_h$ it holds that
	\begin{gather}
	a(u,\widehat{\lambda};v,\widehat\mu) - b(v,\widehat \mu;\phi) = \hat F(v,\widehat \mu), \qquad b(u,\widehat{\lambda};\psi) = 0.
	\end{gather}

The bilinear forms $\widehat a_K$ and $b_K$ are easily seen to satisfy, for all $u^K,v^K \in H^1(K)$, $\hat\mu^K,\hat \lambda^K \in L^2(\partial K)$, $\phi \in L^2(\Sigma)$, the continuity bounds
\begin{equation}\label{contak}
\widehat a_K(u^K,\hat\lambda^K;v^K,\hat\mu^K) \lesssim \left( \| u^K \|_{1,K} + \| \gamma_K^* \hat \lambda^K \|_{-1,K}\right)  \left( \| v^K \|_{1,K} + \| \gamma_K^* \hat \mu^K \|_{-1,K}\right),
\end{equation}
and
\begin{equation}\label{contbk}
b_K(v^K,\hat\mu^K;\phi) \lesssim \| \hat \mu^K \|_{0,\partial K} \| \phi \|_{0,\partial K}.
\end{equation}

\

Setting $\ker b = \{ (v,\hmu) \in \widehat{\mathbb{V}}_h: \ b(v,\hmu;\psi) = 0\ \forall \psi \in \Phi_h \}$, we now observe that, with our choice of the space $\Phi_h$, we have that $(v,\hmu) \in \ker b$ if and only if for some $\mu \in \Lambda_h$ ($\Lambda_h$ as defined in Section \ref{sec:2.2}) it holds  that \begin{equation}\hmuK = \mu(\nu\cdot\nu_K).\label{lambdavshlambda}\end{equation} 
Observe that, for $(u,\hat\lambda),(v,\hat\mu) \in \ker b$, letting $\lambda$ and $\mu$ denote the corresponding elements of $\Lambda_h$ given by \eqref{lambdavshlambda}, it holds that
\begin{equation}\label{equiv}
\hat a(u,\hat \lambda; v, \hat\mu) = a(u,\lambda;v,\mu), \qquad \hat F(v,\hat \mu) = F(v,\mu).
\end{equation}
Then, Lemma \ref{lem:infsup} states an inf-sup condition for  $\widehat a$ on  $\ker b$. Moreover
it is not difficult to prove that 
\begin{equation}\label{2.23}
\inf_{\phi \in \Phi_h} \sup_{(v,\hmu) \in \widehat{\mathbb{V}}_h} \frac{b(v,\hmu;\phi)}{\| \phi \|_{0,\Sigma} \| v,\hmu \|_{\mathbb{V}}} > 0.
\end{equation}
As we are dealing with finite dimensional spaces, for which all norms are equivalent, this is an immediate consequence of the local inf-sup condition
\begin{equation*}
\inf_{{\phi \in \mathbb{P}_\klambda(e)}} \sup_{{\hlambda \in \mathbb{P}_\klambda(e)}} \frac{\int_e \phi\hlambda}{\| \phi \|_{0,e} \| \hlambda \|_{0,e}} \gtrsim 1.
\end{equation*}
As $\widehat{V}_h$ and $\Phi_h$ are finite dimensional spaces, the inf-sup condition for  $\hat a$ on $\ker b$, and the inf-sup condition \eqref{2.23}, together with \eqref{contak} and \eqref{contbk} (as we are in finite dimension these imply continuity with respect to the $\mathbb{V}$ and $L^2(\Sigma)$ norms, though with constants depending on the tessellation), are sufficient to  have existence and uniqueness of the solution to Problem \ref{PbGlobHybrid} (see \cite[Theorem 3.2.1]{Boffi2013}).

\

%
%

It is  easy to realize that Problem  \ref{PbGlob}  and Problem \ref{PbGlobHybrid} are equivalent, and that the solution of the former can be retrieved by actually computing a solution of the latter. Indeed, if $u, \hat \lambda, \phi$ is a solution to Problem \ref{PbGlobHybrid}, then $(u,\hat\lambda) \in \ker b$ and, for the corresponding $\lambda$ given by the relation \eqref{lambdavshlambda}, thanks to \eqref{equiv} it is easy to check that $(u,\lambda)$ is a solution to Problem \eqref{PbGlob}. }

\

{
		It is interesting to give an interpretation of the local stabilization term as the result of a suitable definition of the numerical trace. In the ideal case where $s_K$ is the scalar product for the space $(\Hstar(K))'$ (we recall that
$\Hstar(K) = \{ u \in H^1(K): \ \fint_{\partial K} u = 0\}
$),  
endowed with the norm $| \cdot |_{1,K}$,
it is easy to check that letting $\Riz: (\Hstar(K))' \to \Hstar(K)$ denote the Riesz's isomorphism, which, we recall, is defined in such a way that 
\[
s_K(F,G) = \int_K \nabla \Riz F \cdot \nabla \Riz G = \langle F , \Riz G \rangle = \langle G , \Riz F \rangle,
\]
we have $\Riz = \D^{-1}$.
Considering, for simplicity, the case $t = 0$, we then have
\begin{gather*}
s_K(\D u^K - \gamma_K^* \hat\lambda^K , \gamma_K^* \hat\mu^K )  = \langle \gamma_K^* \hat\mu^K , \Riz(\D u^K - \gamma_K^* \hat\lambda^K) \rangle \\ = \langle  \gamma_K^* \hat\mu^K , u^K - \Riz \gamma_K^* \hat\lambda^K \rangle
=
\int_{\partial K} (u^K - \Riz \gamma_K^* \hat\lambda^K ) \hat\mu^K.
\end{gather*}
The stabilized discrete local problem  (\ref{PbKh1}-\ref{PbKh2}) would then rewrite as
\begin{gather*}
\int_K \nabla u^K \cdot \nabla v^K - \int_{\partial K} \hat\lambda^K v^K = \int_K f v^K,\\
\int_{\partial K}  \widehat u^K \hat\mu^K = \int_{\partial  K}\widehat \phi^K \hat\mu^K
\end{gather*}
with
\[
\widehat u^K   = (1 - \alpha ) u^K + \alpha \gamma_K  \Riz (\gamma_K^* \hat\lambda^K), \qquad \text{ and } \qquad \widehat \phi^K = \phi^K - \alpha \gamma_K  \Riz (f).
\]

\

It is not difficult to check that, if, instead, we define $s_K$ as in \eqref{defskdiscrete}, and we set 
$\tA=(\ta_{ij})$, with $\ta_{ij} = \int_{K} \nabla 
\phi_i \cdot \nabla \phi_j$, then the vector $\vec{x} = (x_i)_{i=1}^N = \tA^{-1} \vec \eta$ would be the vector of coefficient of the function $x = \sum_{i=1 }^N x_i \phi_i \in Y_K$ verifying for all $y \in Y_K$ 
\[
\int_K \nabla x \cdot \nabla y = \int_K \nabla u^K \cdot \nabla y - \int_{\partial K} \hat\lambda^K y.
\]
Considering again the case $t = 0$ and letting $\PiW: \Hstar(K)\to Y_K$ denote the Galerkin projection onto $Y_K$,  we would then be able to rewrite the stabilized problem as
\begin{gather}\label{3.6a}
\int_{K} \nabla u^K \cdot \nabla v^K - \int_{\partial K}
\hat\lambda^K v^K = \int_K f v^K\\
\int_{\partial K}  \widehat u^K \hat\mu^K = \int_{\partial  K}\widehat \phi^K \hat\mu^K, \label{3.6}
\end{gather}
this time with
\[
\widehat u^K  = u^K - \alpha \PiW\Riz(\D u^K - \gamma_K^* \hat\lambda^K ), \qquad \text{ and } \qquad \widehat \phi^K = \phi - \alpha \PiW \Riz f.
\]

Replacing the stiffness matrix $\tA$ with an approximation (as it is done in the Virtual Element Method, when computing $\tA$ as the stiffness matrix relative to the operator $\ta$ defined by \eqref{defta}), results in replacing the Galerkin projection operator $\PiW$ with a spectrally equivalent operator $\widetilde \Pi^K$ 
and setting, in \eqref{3.6},
\begin{equation}\label{defnumtrace}
\widehat u^K  = u^K - \alpha \widetilde \Pi^K \Riz(\D u ^K - \gamma_K^* \hat\lambda^K ), \qquad \text{ and } \qquad \widehat \phi_h^K = \phi - \alpha\widetilde \Pi^K \Riz f.
\end{equation}
More precisely, if $\ta$ is defined by equation \eqref{defta}, it is not difficult to see that, letting $\widetilde \Pi^K: \Hstar(K) \to Y_K$ be defined by
\[
\ta(\widetilde \Pi^K v, y) = \int_K \nabla \Pi^\nabla_K(\widetilde \Pi^K v) \cdot \nabla \Pi^\nabla_K y + \SK(\widetilde \Pi^K v-\Pi^\nabla_K (\widetilde \Pi^K v),y-\Pi^ \nabla_K y) = \int_K \nabla v \cdot\nabla y,
\]
then the local stabilized problem can be rewritten in the form \eqref{3.6a}-\eqref{3.6} with $\widehat u^K$ and $\widehat \phi^K$ given by \eqref{defnumtrace}.

\

{Observe that, unlike what would happen if we used a mesh dipendent stabilization, such as the one proposed in \cite{EwingWangWang} -- that could be interpreted as resulting from defining the numerical trace as a linear combination of the actual trace plus some weighted residual on the fluxes (see \cite{BCMS}) -- the stabilization proposed  here  results in a numerical trace $\widehat u_h^K$ which is indeed the trace of an $H^1(K)$ function.}

}


	\section{Numerical Results}\label{sec:4}

    \begin{figure}
    \centering
    \subfloat{\includegraphics[width=0.3\textwidth]{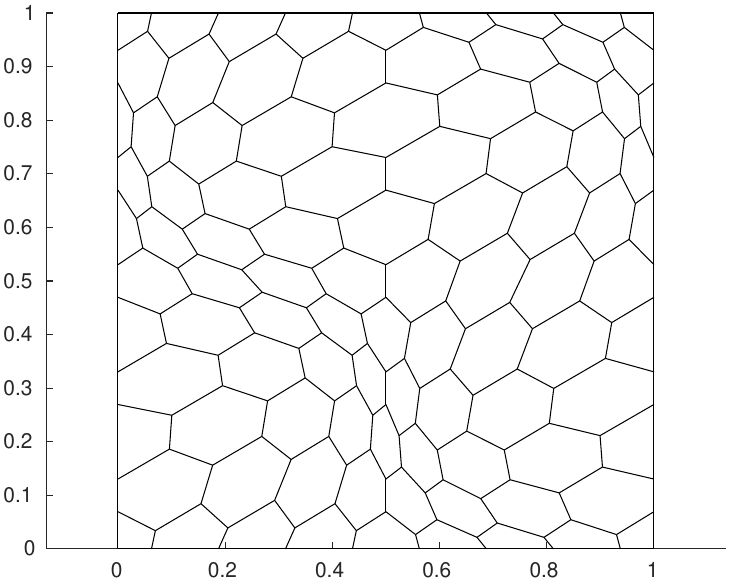}\label{fig:d-hexa}}
    \subfloat{\includegraphics[width=0.3\textwidth]{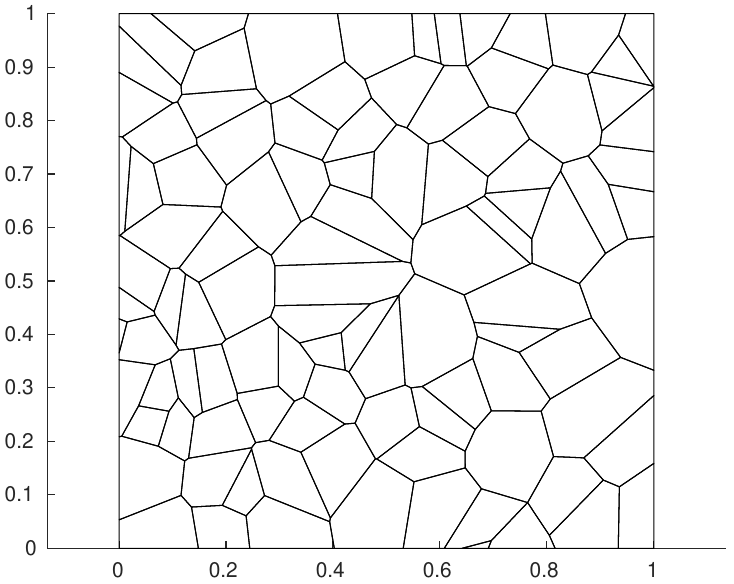}\label{fig:voro}}
    \subfloat{\includegraphics[width=0.3\textwidth]{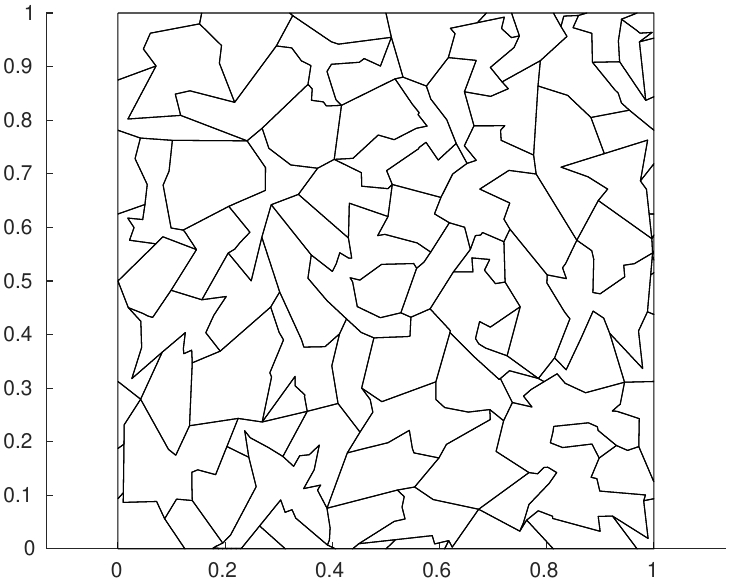}\label{fig:tsp}}
    \caption{(from lest to right) meshes made of deformed hexagons, random Voronoi cells, and random polygons.}
    \end{figure}

	\begin{table}
\centering
\caption{Meshes of deformed hexagons.}
\label{tab:mesh-d-hexa}
\begin{tabular}{
c
S[table-format=4.0]
S[table-format=5.0]
S[table-format=1.{\roundPrecision}e-1]
S[table-format=1.{\roundPrecision}e-1]
S[table-format=1.{\roundPrecision}e-1]
S[table-format=1.{\roundPrecision}e+1]
S[table-format=1.2]
S[table-format=2.0]
}
\toprule
Mesh & {$N_\text{el}$} & {$N_\text{ed}$} & {$h_\text{max}$} & {$h_\textup{min}$} &{$h_\text{av}$}
&${\gamma_0}$ &{$\gamma_1$}& {$N^\star$}\\
\midrule
d-hexa$_{1}$ & 822 & 2467 & 7.113742e-02 & 8.794092e-03 & 2.01e-02
& 1.209727e+01 & 5.052029 & 6\\
d-hexa$_{2}$ & 1415 & 4246 & 5.424294e-02 & 6.819149e-03 & 1.53e-02
& 1.225613e+01 & 4.934710& 6\\
d-hexa$_{3}$ & 2270 & 6811 & 4.233760e-02 & 5.124184e-03 & 1.21e-02
& 1.240315e+01 & 5.163088& 6\\
d-hexa$_{4}$ & 3203 & 9610 & 3.566680e-02 & 4.397059e-03 & 1.02e-02
& 1.223047e+01 & 5.052014& 6\\
d-hexa$_{5}$ & 4296 & 12889 & 3.082316e-02 & 3.827593e-03 & 8.81e-03
& 1.230091e+01 & 5.002315& 6\\
d-hexa$_{6}$ & 5711 & 17134 & 2.659539e-02 & 3.232217e-03 & 7.64e-03
& 1.247704e+01 & 5.130111& 6\\
\bottomrule
\end{tabular}
\end{table}

\begin{table}
\centering
\caption{Meshes of random Voronoi cells.}
\label{tab:mesh-voro}
\begin{tabular}{
c
S[table-format=5.0]
S[table-format=6.0]
S[table-format=1.{\roundPrecision}e-1]
S[table-format=1.{\roundPrecision}e-1]
S[table-format=1.{\roundPrecision}e-1]
S[table-format=1.{\roundPrecision}e+1]
S[table-format=1.{\roundPrecision}e+1]
S[table-format=2.0]
}
\toprule
Mesh & {$N_\text{el}$} & {$N_\text{ed}$} & {$h_\text{max}$} & {$h_\textup{min}$} &{$h_\text{av}$}&${\gamma_0}$& {$\gamma_1$}&{$N^\star$}\\
\midrule
voro$_{1}$ & 2500 & 7505 & 6.384666e-02 & 6.192528e-06 &1.15e-02
& 1.422215e+01 & 5.852731e+03&11\\
voro$_{2}$ & 5000 & 15007 & 4.344562e-02 & 5.845179e-07 & 8.17e-03
& 1.456663e+01 & 3.404522e+04&14\\
voro$_{3}$ & 10000 & 30006 & 3.470002e-02 & 1.732139e-07 & 5.77e-03
& 2.525383e+01 & 9.493668e+04&12\\
voro$_{4}$ & 20000 & 60010 & 2.405393e-02 & 2.138871e-07 &4.08e-03
& 2.087832e+01 & 7.246942e+04&13\\
voro$_{5}$ & 40000 & 120006 & 1.726980e-02 & 8.256465e-08 &2.89e-03
& 2.675024e+01 & 7.178744e+04&13\\
voro$_{6}$ & 80000 & 240027 & 1.140086e-02 & 5.998477e-09 &2.04e-03
& 2.881822e+01 & 1.100228e+06&13\\
\bottomrule
\end{tabular}
\end{table}

\begin{table}
\centering
\caption{Meshes of random polygons. For this class of meshes, $\gamma_0$ is not computed.}
\label{tab:mesh-tsp}
\begin{tabular}{
c
S[table-format=5.0]
S[table-format=6.0]
S[table-format=1.{\roundPrecision}e-1]
S[table-format=1.{\roundPrecision}e-1]
S[table-format=1.{\roundPrecision}e-1]
S[table-format=1.{\roundPrecision}e+1]
S[table-format=2.0]
}
\toprule
Mesh & {$N_\text{el}$} & {$N_\text{ed}$} & {$h_\text{max}$} & {$h_\textup{min}$} &{$h_\text{av}$}& {$\gamma_1$}&{$N^\star$}\\
\midrule
tsp$_{1}$ & 1000 & 7367 & 1.206953e-01 & 8.610327e-05 &1.83e-02
& 5.607546e+02& 36\\
tsp$_{2}$ & 2000 & 14197 & 7.574066e-02 & 4.323479e-05 &1.29e-02
& 1.335083e+03& 41\\
tsp$_{3}$ & 4000 & 28003 & 5.090225e-02 & 1.338424e-05 &9.13e-03
& 1.614139e+03& 41\\
tsp$_{4}$ & 8000 & 54668 & 4.118779e-02 & 7.039107e-06 & 6.46e-03
& 8.729327e+02& 40\\
tsp$_{5}$ & 16000 & 107550 & 2.654308e-02 & 1.681474e-06 & 4.56e-03 
& 3.542747e+03& 46\\
tsp$_{6}$ & 32000 & 213570 & 1.998982e-02 & 5.452284e-07 & 3.23e-03
& 1.564320e+03& 65\\
\bottomrule
\end{tabular}
\end{table}

	We take the domain $\Omega$ to be the unit square $[0,1]\times[0,1]$. 
	We solve Problem  \ref{Pb:strong} with Dirichlet boundary data $g$, and load term $f$ chosen in such a way that
	\[
	u = \frac{1}{128\pi^2}\cos(8\pi x)\cos(8\pi y)
	\]
	is the exact solution. The stabilization parameters are chosen to be $\alpha = t = 1$. {We test our method on three sequences of meshes with increasingly degrading shape regularity}: deformed hexagonal meshes (test case 1, Figure~\ref{fig:d-hexa}), random Voronoi meshes (test case 2, Figure~\ref{fig:voro}), and meshes made of random polygons (test case 3, Figure~\ref{fig:tsp}) generated as follows: i) throw random points inside $\Omega$; ii) partition them into a given number of clusters; iii) join the points of each cluster with the shortest closed tour, i.e., solve the Traveling Salesman Problem; iv) mesh the complement of the polygons obtained at step iii) with triangles and agglomerate them. Geometrical data for these meshes are shown in Tables~\ref{tab:mesh-d-hexa}, \ref{tab:mesh-voro}, \ref{tab:mesh-tsp}, respectively. For each mesh, we provide: $N_\text{el}$, the number of elements of $\mathcal T_h$; $N_\text{ed}$, the number of edges of $\mathcal T_h$; $h_\text{max} = \max_{K\in\Omega_h} h_K$, the maximum element diameter; $h_\textup{min} = \min_{K\in\mathcal T_h} h_{\textup{min},K}$, where $h_{\textup{min},K}$ is the minimum distance between any two vertices of $K$; 
	$h_\text{av} = N_\text{el}^{-1/2}$, an estimate of the average mesh-size;
	$\gamma_0 = \max_{K\in\mathcal T_h}\frac{h_K}{\rho_K}$, where $\rho_K$ is the radius of the largest circle that is contained inside $K$; $\gamma_1 = \max_{K\in\mathcal T_h}\frac{h_K}{h_{\textup{min},K}}$; $\Nstar = \max_K N_K$, the maximum of the number of edges in each element.

\

   In order to compute $u$ and $\lambda$ we solve the equivalent hybridized Problem \ref{PbGlobHybrid}. Since, for each $K$, \eqref{PbK1}-\eqref{PbK2} yield a local discrete Dirichlet problem, we can resort to a {static condensation} procedure, allowing to reduce the size of the resulting algebraic equation by expressing $u^K, \hat \lambda^K$ as a function of the sole variable $\varphi|_{\partial K}$. At this point, we can use \eqref{Pbcoupling}, which imposes continuity of the fluxes $\lambda$, to glue all the local problems together and obtain a global system of equations where only $\varphi$ appears as unknown.    
     In the present tests the global system is solved with the direct solver STRUMPACK~\cite{strumpack}. Reconstruction of $u, \hat\lambda$ is done by solving local problems in parallel. Remark that other approaches yielding efficient implementation can be consideres (see, for instance \cite{Arayaetal}).

\

\begin{figure}
   	\includegraphics[width=8.25cm]{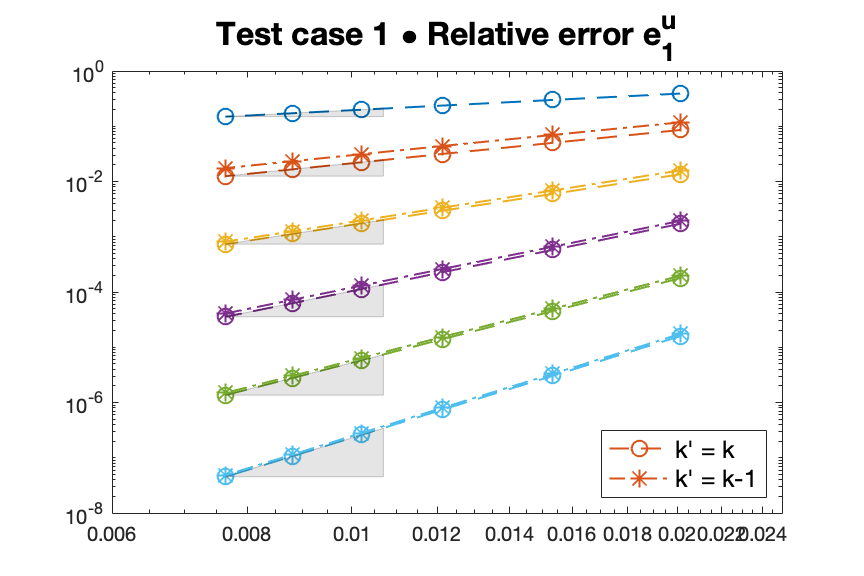}
   	\includegraphics[width=8.25cm]{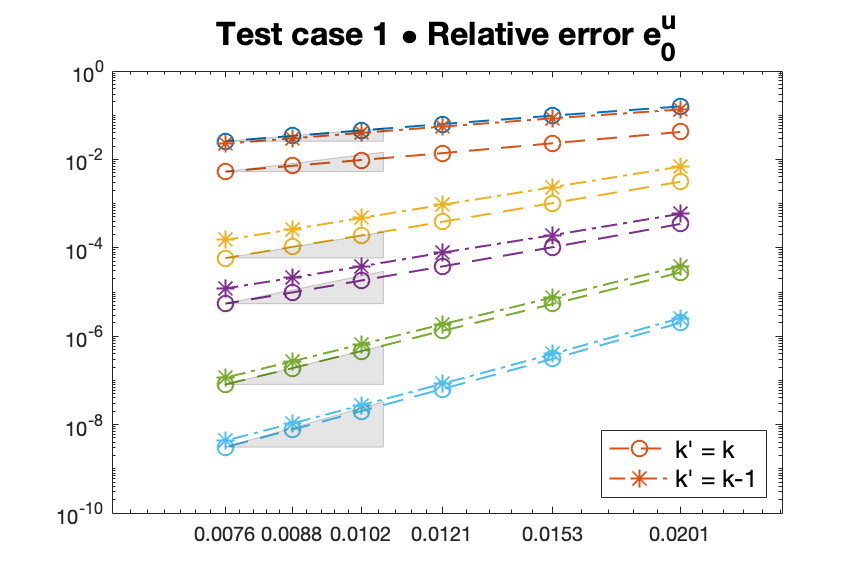}
   \caption{Test Case 1. $\| \cdot \|_{1,*}$ (left) and $\| \cdot \|_{0,\Omega}$ (right) errors for $k=\klambda=1,\cdots,6$ and $k = \klambda + 1 = 2,\cdots,6$.  }
  \label{figtest1}
   \end{figure}

\begin{figure}
	\includegraphics[width=8.25cm]{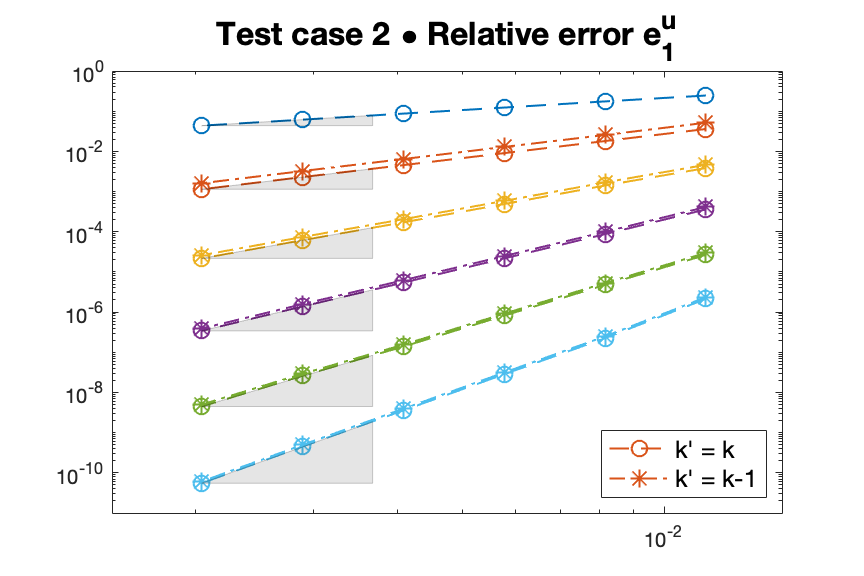}
	\includegraphics[width=8.25cm]{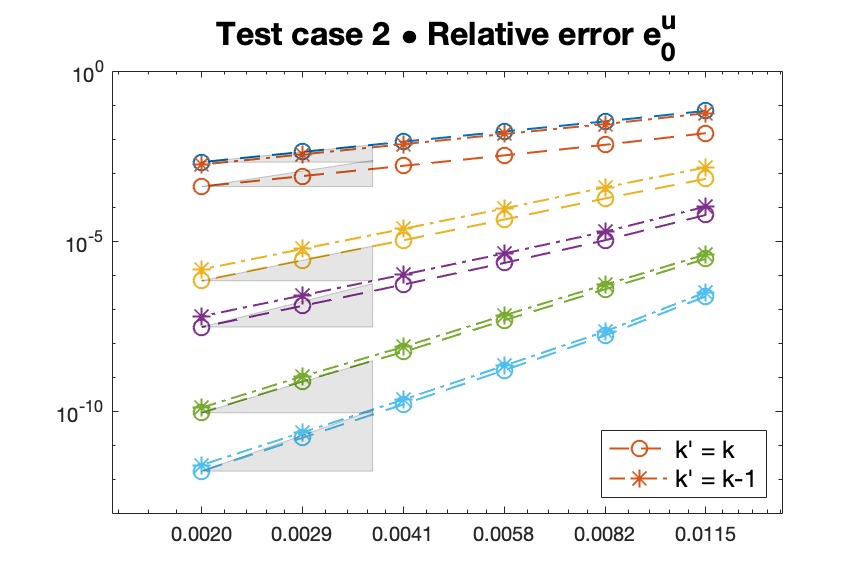}
	\caption{Test Case 2. $\| \cdot \|_{1,*}$ (left) and $\| \cdot \|_{0,\Omega}$ (right) errors for $k=\klambda=1,\cdots,6$ and $k = \klambda + 1 = 2,\cdots,6$.  }  \label{figtest2}
\end{figure}

\begin{figure}
	\includegraphics[width=8.25cm]{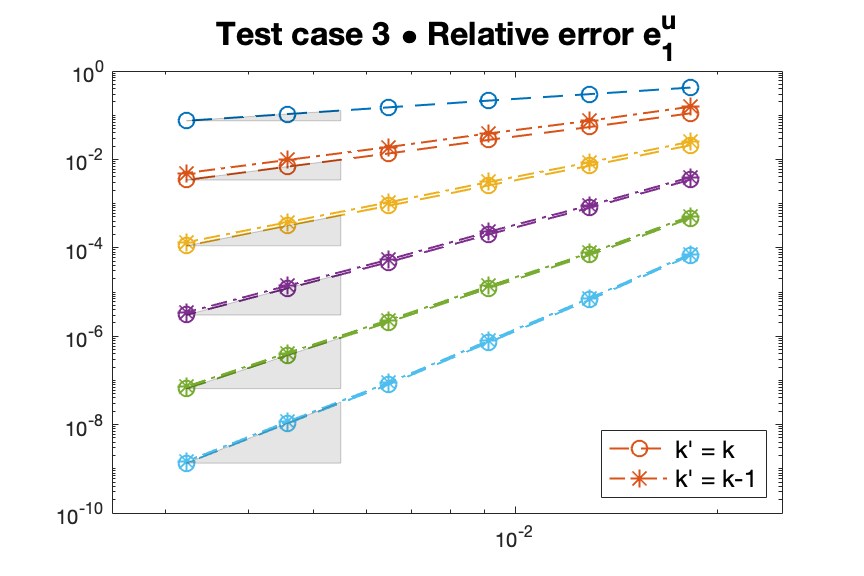}
	\includegraphics[width=8.25cm]{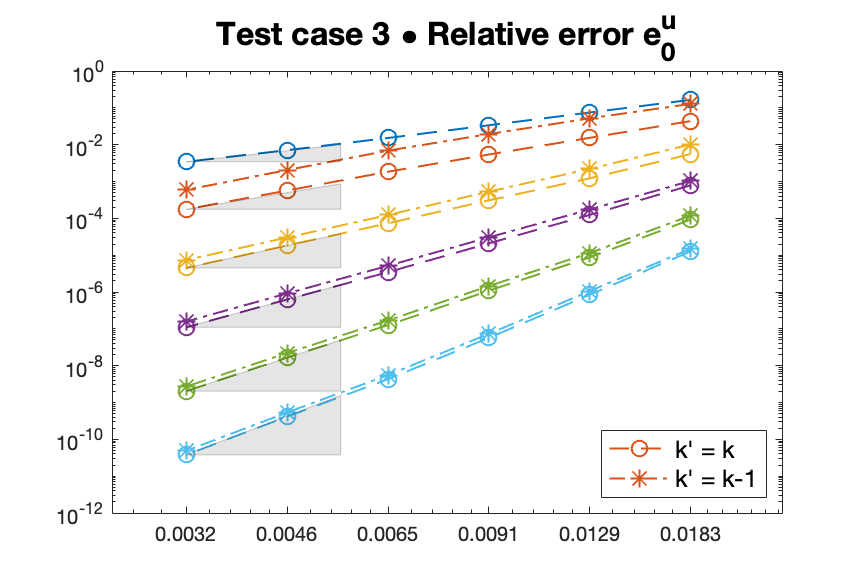}
	\caption{Test Case 3. $\| \cdot \|_{1,*}$ (left) and $\| \cdot \|_{0,\Omega}$ (right) errors for $k=\klambda=1,\cdots,6$ and $k = \klambda + 1 = 2,\cdots,6$.  }  \label{figtest3}
\end{figure}

    For the three test cases, Figures \ref{figtest1}, \ref{figtest2} and \ref{figtest3} respectively show the relative errors 
    \[e^u_1 = \frac{\| u - u_h \|_{0,\Omega} + | u - u_h |_{1,*}}{\| u \|_{0,\Omega} + | u |_{1,\Omega}}, \qquad e^u_0 = \frac{\|u-u_h \|_{0,\Omega}}{\| u \|_{0,\Omega}} \]
   (where $| \cdot |_{1,*}$ denotes the broken $H^1$ seminorm), for $k=\klambda=1,\cdots,6$ (dotted lines with circular markers) and $k = \klambda+1 = 2,\cdots,6$ (dashed line with asterisks markers)
(for technical reasons, related to the  actual implementation of the stabilization term, we did not test the case $k=1$, $\klambda = 0$).   
The errors are plotted, in logarithmic scale, against $h_{\text{av}} = 1/\sqrt{N_\text{el}}$ (which ideally behaves as an average element size). The slope of the gray triangles in the pictures shows the optimal convergence rate attainable by the best approximation (equals to $k$ for $e_1^u$ and $k+1$ for $e_0^u$).

	\
	
We observe that the results confirm the theoretical estimate, with the correct order of convergence for the broken $H^1$ norm of the error, i.e. $\mathcal O(h^k)$, as $h$ tends to zero. Observe also that, as far as the choice of $\klambda$ is concerned, when considering the $H^1$ norm, there is very little difference between $\klambda = k$ and $\klambda = k-1$.

 As far as the convergence in the $L^2$ norm is concerned, observe that, for the first two  test cases we get the optimal convergence rate only for the even values of $k$. This is consistent with results obtained for non symmetric interior penalty approximations of linear elliptic problems~\cite{houston2002,BABUSKA1999103}. Surprisingly, for the third test case it appears that the method gets very close to  the optimal order of convergence also for the odd values of $k$. Also to be remarked is the fact that, when considering the error in the $L^2$ norm, the discretization with $\klambda=k$ behaves sensibly better than the one with $\klambda=k-1$, at least for low values of $k$. In particular, for the first two test cases, the curve relative to the discretozation with $k=k'=1$ is superposed to the one relatve to $k = 2$, $k'=1$. If, for these two cases, we compare the number of degrees of freedom, we realize that the discretization with $\klambda = k = 1$ allows to attain the same $L^2$ error as the one with $k=2$, $k'=1$, with $4 \times N_\text{el}$ less degrees of freedom. In general, for a fixed $k$ using $\klambda = k$ yields an error $3$ to $4$ times smaller than the one obtained with $\klambda = k-1$, with $1$ extra degree of freedom per edge.

	\section{Conclusions}
We presented and analyzed a hybrid Discontinuous Galerkin method on a polygonal tessellation for the Poisson problem in two dimensions, with a new design for the stabilization term, based on an algebraic representation for the scalar product of the duals of the spaces $H^1(K)$ for $K$ element of the tessellation. Following the general recipe provided in \cite{AlgStab}, such scalar products can, in fact, be numerically realized, via the introduction of a (minimal) auxiliary space, for which no approximation properties are required but which has to satisfy an inf-sup condition.
Under quite weak shape regularity assumptions, allowing for the presence of elements with very small edges, we proved optimal error estimates (confirmed by the results of the numerical tests), thus demonstrating the feasibility and the potential of a stabilization approach where some residual term is penalized in the norm of the dual space where it naturally lives. We believe that such an approach can potentially be applied to replace the mesh dependent stabilization terms appearing also in  other discontinuous Galerkin formulations, and, possibly, also beyond the framework of discontinuous Galerkin methods.

	\bibliographystyle{plain}
	\bibliography{biblio}

\end{document}